\documentclass [10pt,a4paper] {article}
\usepackage{color}
\usepackage{array}
\usepackage{tabularx}
\usepackage{amsmath}
\usepackage{amsfonts}
\usepackage{bbm}
\usepackage{latexsym}
\usepackage{epsfig}
\newtheorem{remark}{Remark}
\newtheorem{lemma}{Lemma}
\newtheorem{theorem}{Theorem}
\newtheorem{definition}{Definition}
\newtheorem{proposition}{Proposition}
\newtheorem{corollary}{Corollary}
\newtheorem{example}{Example}

\begin{document}

\setcounter{page}{25}

\newcommand{\proof}{\paragraph{Proof}}
\newcommand{\fin}{$\Box$\\}
\newcommand{\ds}{\displaystyle}
\newcommand{\saut}[1]{\hfill\\[#1]}
\newcommand{\vsp}{\vspace{.15cm}}
\newcommand{\difrac}{\displaystyle \frac}
\newcommand{\dist}{\textrm{dist}}
\newcommand{\mbf}{\textbf}
\title{From $N$ parameter fractional Brownian motions to 
$N$ parameter multifractional Brownian motions}
\author{Erick Herbin \vspace{5mm} \\ 
\small INRIA,
Domaine de Voluceau, Rocquencourt,
BP 105,
78153 Le Chesnay Cedex, France \\
\small erick.herbin@inria.fr\\
\small and\\
\small Dassault Aviation, 
78 quai Marcel Dassault,
92552 Saint-Cloud Cedex, France \\
\small erick.herbin@dassault-aviation.fr}
\date{November 2002}

\maketitle
\begin{abstract}
Multifractional Brownian motion is an extension of the well-known fractional
Brownian motion where the H\"older regularity is allowed to vary along the
paths. In this paper, two kind of multi-parameter extensions of mBm are studied:
one is isotropic while the
other is not. For each of these processes, a moving average representation, 
a harmonizable representation, and the covariance structure are given.\\
The H\"older regularity is then studied. In particular, the case of an irregular
exponent function $H$ is investigated.
In this situation,
the almost sure pointwise and local H\"older exponents of the multi-parameter
mBm are proved to be equal to the correspondent exponents of $H$.
Eventually, a local asymptotic self-similarity property is proved. The limit
process can be another process than fBm.
\end{abstract}

{\sl AMS classification\/}: 62\,G\,05, 60\,G\,15, 60\,G\,17, 60\,G\,18.

{\sl Keywords\/}: fractional Brownian motion, Gaussian processes, H\"older
regularity, local asymptotic self-similarity, multi-parameter processes.


\section{Introduction}
In many applications, fractional Brownian motion (fBm) seems to fit very well
to random phenomena. 
Recall that it can be defined by one of the four following properties.
Let $H\in (0,1)$ ($H$ is sometimes called the Hurst parameter).
\begin{itemize}
\item $B^H$ is a centered Gaussian process such that
\[
\forall s,t\in\textbf{R}_{+};\;
E\left[B^H_s B^H_t\right]=\frac{1}{2}\left[s^{2H}+t^{2H}-|t-s|^{2H}\right]
\]
\item the process $B^H$ such that
\[
\forall t\in\textbf{R}_{+};\;
B^H_t=\int_{-\infty}^0 \left[(t-u)^{H-\frac{1}{2}}-(-u)^{H-\frac{1}{2}}\right].\mathbb{W}(du)
+\int_0^t (t-u)^{H-\frac{1}{2}}.\mathbb{W}(du)
\]
is a fBm,
\item the process $B^H$ such that
\[
\forall t\in\textbf{R}_{+};\;
B^H_t=\int_{\textbf{R}}\frac{e^{it\xi}-1}{\left|\xi\right|^{H+\frac{1}{2}}}.\hat{\mathbb{W}}(d\xi)
\]
is a fBm,
\item $B^H$ is the unique self-similar Gaussian process with stationary
increments.
\end{itemize}

Its efficiency has already been shown in simulation of 
traffic on Internet or in finance. This induced some recent progress such as
stochastic integration against fBm.\\
However, the main limitation of fBm is that the H\"older regularity is constant
along the paths. \\
Multifractional Brownian motion (mBm) has been independently introduced in
\cite{benassi} and \cite{peltier}. This process is a generalization
of fractional Brownian motion where the Hurst parameter $H$
is substituted by a function $t\mapsto H(t)$. As a consequence the H\"older 
exponent is allowed to vary along trajectories.\\
The different definitions by the two groups of authors provided two different 
representations of mBm.

Peltier and Levy-Vehel (\cite{peltier}) defined the mBm from the moving average definition
of the fractional Brownian motion 
\[
X_t=\int_{-\infty}^0 \left[(t-u)^{H(t)-\frac{1}{2}}-(-u)^{H(t)-\frac{1}{2}}\right].\mathbb{W}(du)
+\int_0^t (t-u)^{H(t)-\frac{1}{2}}.\mathbb{W}(du)
\]
where $t\mapsto H(t)$ is a H\"older function.

Benassi, Jaffard and Roux (\cite{benassi}) defined the mBm from the harmonizable representation
of the fBm
\[
X_t=\int_{\textbf{R}}\frac{e^{it\xi}-1}{\left|\xi\right|^{H(t)+\frac{1}{2}}}.\hat{\mathbb{W}}(d\xi)
\]

These two definitions were proved to be equivalent up to a multiplicative
deterministic function (\cite{cohencov}).

Moreover, in \cite{cohen} the covariance function of this Gaussian process has been proved
to be 
\[
E\left[X_s X_t\right]=D\left(H(s),H(t)\right)\left[|s|^{H(s)+H(t)}+|t|^{H(s)+H(t)}
-|t-s|^{H(s)+H(t)}\right]
\]
where $D$ is a known deterministic function.\\
The goal of this paper is to study some multi-parameter extension of the 
multifractional Brownian motion, ie a stochastic process indexed by 
$\textbf{R}_{+}^N$, which is an mBm when $N=1$. One extension has already been
considered in \cite{benassi}.\\
2D extension of fractional Brownian motion has been already used in various
applications such as underwater terrain modeling (\cite{beatrice}). 
It may be more realistic to allow local 
regularity to vary at each point : our extension of mBm in $\textbf{R}^2$ may be 
used for this kind of application.

\section{Multi-parameter extension of the fractional Brownian motion}
Since multifractional Brownian motion is an extension of fractional Brownian
motion, we start with a review of the existing extensions of fBm.
Most of the results in this section are well-known, but we give new proofs
based only on the covariance functions.

In the same way as Brownian motion has two main multi-parameter extensions: 
Levy Brownian motion and Brownian sheet,
two different multi-parameter extensions of fractional Brownian motion
have been defined.

\subsection{Levy fractional Brownian motion}
This process can be seen as an isotropic extension of the fractional Brownian 
motion.
Indeed, for the fBm, we have for all $s,t\in\textbf{R}_{+}$
\[
E\left[X_t - X_s\right]^2=|t-s|^{2H}
\]

A natural idea to extend this process for a set of index $\mathcal{T}\subset\textbf{R}_{+}^N$
is to substitute the absolute value by a norm.
We get the Levy fractional Brownian motion, which is defined to be a centered
Gaussian process of covariance function
\begin{equation} \label{levycov}
E\left[X_s X_t\right]=\frac{1}{2}\left[\|s\|^{2H}+\|t\|^{2H}-\|t-s\|^{2H}\right]
\end{equation}
 
There are several definitions of this process by its trajectories.
Among these, it can be defined as integral against white noise.
Lindstrom stated the following (see \cite{lindstr}).

\begin{proposition}\label{levymoyrep}
The process defined by
\begin{eqnarray}\label{levymob}
X_t=\int_{\textbf{R}^N}\left[\|t-u\|^{H-\frac{N}{2}}-\|u\|^{H-\frac{N}{2}}\right]
\mathbb{W}(du)
\end{eqnarray}
is a Levy fractional Brownian motion up to a multiplicative constant.
\end{proposition}
\proof
This process is obviously Gaussian and centered. Thus we only have to show that
the covariance function is of the form (\ref{levycov}).
We have
\begin{eqnarray*}
E\left[X_s - X_t\right]^2&=&\int_{\textbf{R}^N}\left[\|t-u\|^{H-\frac{N}{2}}-\|s-u\|^{H-\frac{N}{2}}\right]^2.du\\
&=&\int_{\textbf{R}^N}\left[\|t-s-u\|^{H-\frac{N}{2}}-\|u\|^{H-\frac{N}{2}}\right]^2.du
\end{eqnarray*}
We consider the change of variables from $\textbf{R}^N$ into itself, 
$v=\phi(u)$,
where $\phi$ is the linear application which maps the canonic basis of $\textbf{R}^N$
to the orthonormal basis $\left(e_1=\frac{t-s}{\|t-s\|},e_2,\dots,e_N\right)$.
The differential of $\phi$ in any $u\in\textbf{R}^N$ is itself and the Jacobian
\begin{eqnarray*}
J\phi^{-1}_{(v_1,\dots,v_N)}=\left|\det \left(\phi^{-1}\right)\right|
=1
\end{eqnarray*}
because the matrix of $\phi$ is orthogonal.\\
We have
\begin{eqnarray*}
\|t-s-u\|^2&=&(t-s-u)^2=\|t-s\|^2-2<t-s,u>+\|u\|^2\\
&=&\|t-s\|^2-2\|t-s\|.v_1+\|v\|^2\\
&=&\left(\|t-s\|.\epsilon_1-v\right)^2
\end{eqnarray*}
We obtain
\begin{eqnarray*}
E\left[X_s - X_t\right]^2&=&\int_{\textbf{R}^N}\left[\| \; \|t-s\|.\epsilon_1-v\|^{H-\frac{N}{2}}-\|v\|^{H-\frac{N}{2}}\right]^2.dv\\
&=&\|t-s\|^{2H-N}\int_{\textbf{R}^N}\left[\|\epsilon_1-\frac{v}{\|t-s\|}\|^{H-\frac{N}{2}}-\|\frac{v}{\|t-s\|}\|^{H-\frac{N}{2}}\right]^2.dv
\end{eqnarray*}
and after the second change of variables, 
\[
v=\|t-s\|.w=\|t-s\| Id.w
\]
we get
\begin{eqnarray*}
E\left[X_s - X_t\right]^2=\|t-s\|^{2H}\underbrace{\int_{\textbf{R}^N}\left[\|\epsilon_1-w\|^{H-\frac{N}{2}}-\|w\|^{H-\frac{N}{2}}\right]^2.dw}
_{K_{N,H}}
\end{eqnarray*}
therefore
\begin{eqnarray*}
E\left[X_s X_t\right]=K_{N,H}\left[\|s\|^{2H}+\|t\|^{2H}-\|t-s\|^{2H}\right]
\end{eqnarray*}

\fin

The harmonizable representation of fractional Brownian motion can also be
generalized. Before that, let's recall briefly definitions of white noise
and its Fourier transform.\\
In the following, we will denote $L^2_{\textbf{C}}(\textbf{R}^N)$ the
set of functions $f:\textbf{R}^N\rightarrow\textbf{C}$ such that
$\int_{\textbf{R}^N}|f(u)|^2 du <\infty$.
\begin{definition}
The complex isonormal process is defined to be a centered Gaussian process
$W=\left\{W(f); f\in L^2_{\textbf{C}}(\textbf{R}^N)\right\}$ such that
\[
\forall f,g\in L^2_{\textbf{C}}(\textbf{R}^N);\;
E\left[W(f)\overline{W(g)}\right]=\int_{\textbf{R}^N} f(u)\overline{g(u)}.du
\]
Then, white noise $\mathbb{W}$ can be defined by
\[
\mathbb{W}(E)=W\left(\mathbbm{1}_E\right)
\]
\end{definition}

\begin{definition}
A Gaussian process $\left\{\hat{W}(f);f\in L^2_{\textbf{C}}(\textbf{R}^N)\right\}$ 
is said to be the Fourier transform of a complex
isonormal process $\left\{W(f);f\in L^2_{\textbf{C}}(\textbf{R}^N)\right\}$ 
if for all $f\in L^2_{\textbf{C}}(\textbf{R}^N)$
\[
\hat{W}(f)=W(\hat{f})
\]
where $\hat{f}$ if the Fourier transform of the function $f$.\\
The Fourier transform of white noise is defined in the same way.

\end{definition}

This complex measure is usually used to define
the harmonizable representation of fractional Brownian motion
\[
B^H_t=\int_{\textbf{R}}\frac{e^{it\xi}-1}{\left|\xi\right|^{H+\frac{1}{2}}}.\hat{\mathbb{W}}(d\xi)
\]
that can be generalized in the following.\\

\begin{proposition}
The process defined by
\begin{equation}
X_t=\int_{\textbf{R}^N}\frac{e^{i\left<t,\xi\right>}-1}{\|\xi\|^{H+\frac{N}{2}}}.\hat{\mathbb{W}}(d\xi)
\end{equation}
where $\hat{\mathbb{W}}$ is the Fourier transform of white noise in $\textbf{R}^N$,\\
is a Levy fractional Brownian motion up to a multiplicative constant.
\end{proposition}
\proof
As will be done for multifractional
Brownian field, the Fourier transform of the kernel of representation
(\ref{levymob}) could be directly computed.
But as this representation defines a real centered Gaussian process,
it is enough to show that the covariance function has the form (\ref{levycov}).\\
For all $t\in\textbf{R}^N$, let's denote by $f_t$ the function
$\xi\mapsto\frac{e^{i<t,\xi>}-1}{\|\xi\|^{H+\frac{N}{2}}}$
and consider the centered Gaussian process 
$X=\left\{X_t=\hat{W}(f_t); t\in\textbf{R}_{+}^N\right\}$.\\
First of all, let's show that, almost surely, $\hat{W}(f_t)\in\textbf{R}$.\\
In fact, using $\hat{W}(f_t)=W(\hat{f_t})$, showing that
$\hat{f_t}\in\textbf{R}$ is sufficient.\\
Indeed, by 
$E\left[Im\left(W(f)\right)\right]^2=\int_{\textbf{R}^N}\left(Im\left(f\right)\right)^2=0$,
$f\in\textbf{R}$ imply $W(f)\in\textbf{R}$ almost surely, and
\begin{eqnarray*}
Im(\hat{f_t}(x))&=&\int_{\textbf{R}^N}Im\left(e^{-i<x,u>}f_t(u)\right).du\\
&=&\int_{\textbf{R}^N}\frac{\sin<t-x,u>}{\|u\|^{H+\frac{N}{2}}}.du\\
&=&0
\end{eqnarray*}
by parity.\\
The process $X$ is therefore real and its covariance function is
\begin{eqnarray*}
E\left[X_s X_t\right]&=&E\left[\hat{W}(f_s)\overline{\hat{W}(f_t)}\right]\\
&=&\int_{\textbf{R}^N}\frac{\left(e^{i<s,\xi>}-1\right) 
\left(e^{-i<t,\xi>}-1\right)}{\|\xi\|^{2H+N}}.d\xi\\
&=&\int_{\textbf{R}^N}\frac{e^{i<s-t,\xi>}-e^{i<s,\xi>}-e^{-i<t,\xi>}+1}{\|\xi\|^{2H+N}}.d\xi
\end{eqnarray*}
Then we have to consider $3$ integrals of the form 
$\int_{\textbf{R}^N}\frac{1-e^{i<t,\xi>}}{\|\xi\|^{2H+N}}.d\xi$.\\
As in proposition \ref{levymoyrep},
for $t\in\textbf{R}^N$ fixed, consider the change of variables from $\textbf{R}^N$
into itself, $u=\phi\left(\xi\right)$
where $\phi$ is the linear application which maps the canonic basis of $\textbf{R}^N$
to the orthonormal basis $\left(e_1=\frac{t}{\|t\|},e_2,\dots,e_N\right)$.\\
Then, we get
\[
\int_{\textbf{R}^N}\frac{1-e^{i<t,\xi>}}{\|\xi\|^{2H+N}}.d\xi
=\int_{\textbf{R}^N}\frac{1-e^{i\|t\|.u_1}}{\|u\|^{2H+N}}.du
\]
After the second change of variables
\begin{eqnarray*}
v=\|t\|.u=\|t\|Id.u\\
dv=\|t\|^N.du
\end{eqnarray*}
we get
\begin{eqnarray*}
\int_{\textbf{R}^N}\frac{1-e^{i<t,\xi>}}{\|\xi\|^{2H+N}}.d\xi
=\frac{\|t\|^{2H+N}}{\|t\|^{N}}\underbrace{\int_{\textbf{R}^N}\frac{1-e^{iv_1}}{\|v\|^{2H+N}}.dv}_{C_{N,H}>0}
\end{eqnarray*}
Proceeding the same way for the $2$ other integrals, we can conclude
\begin{eqnarray*}
E\left[X_s X_t\right]=C_{N,H}\left[\|s\|^{2H}+\|t\|^{2H}-\|t-s\|^{2H}\right]
\end{eqnarray*}
which shows that the process 
$\left\{\frac{1}{\sqrt{C_{N,H}}}\hat{W}(f_t), t\in\textbf{R}_{+}^N\right\}$
is a Levy fractional Brownian motion.
\fin

\subsection{Fractional Brownian sheet}
On the contrary to the Levy fractional Brownian motion, this process is not
isotropic. In particular, we can have different Hurst parameters in each
of the $N$ directions.\\
For the fBm, we have for all $s,t\in\textbf{R}_{+}$
\[
E\left[X_s X_t\right]=\frac{1}{2}\left[s^{2H}+t^{2H}-|t-s|^{2H}\right]
\]
As in the definition of Brownian sheet,
another way to generalize fBm is to set the covariance equal to the tensor 
product of one dimensional covariances. 
Then, fractional Brownian sheet (fBs) is defined to be a centered Gaussian
process of covariance function
\begin{equation} \label{fbs}
E\left[X_s X_t\right]=\prod_{i=1}^{N}\frac{1}{2}\left(s_i^{2H_i}+t_i^{2H_i}
-|t_i-s_i|^{2H_i}\right)
\end{equation}

As in the isotropic case, this process has two different representations 
by its trajectories.

\begin{proposition} \label{fbsmm}
The process defined by
\[
X_t=\int_{\textbf{R}^N} \prod_{i=1}^N \left[|t_i-u_i|^{H_i-\frac{1}{2}}-|u_i|^{H_i-\frac{1}{2}}\right]\mathbb{W}(du)
\]
is a fractional Brownian sheet, up to a multiplicative constant.
\end{proposition}

\begin{remark}
In \cite{legerpontier}, Pontier/Leger introduced another
moving average representation of fractional Brownian sheet.
\[
X_t=\int_{\textbf{R}^N} \prod_{i=1}^N \left[(t_i-u_i)_{+}^{H_i-\frac{1}{2}}-(-u_i)_{+}^{H_i-\frac{1}{2}}\right]\mathbb{W}(du)
\]
\end{remark}

\proof
This process is obviously Gaussian and centered. Thus, we only need to show that
its covariance function has the expected form.
We compute
\begin{eqnarray*}
E\left[X_s X_t\right]&=&\int_{\textbf{R}^N} \prod_{i=1}^N \left[|s_i-u_i|^{H_i-\frac{1}{2}}-|u_i|^{H_i-\frac{1}{2}}\right]
\left[|t_i-u_i|^{H_i-\frac{1}{2}}-|u_i|^{H_i-\frac{1}{2}}\right].du\\
&=&\prod_{i=1}^N\int_{\textbf{R}}\left[|s_i-u_i|^{H_i-\frac{1}{2}}-|u_i|^{H_i-\frac{1}{2}}\right]
\left[|t_i-u_i|^{H_i-\frac{1}{2}}-|u_i|^{H_i-\frac{1}{2}}\right].du_i
\end{eqnarray*}
We can see that the factor corresponding to each $i$, is the covariance of a fBm 
with Hurst parameter $H_i$ (or a Levy fractional Brownian motion with $N=1$).
Then we have
\begin{eqnarray*}
E\left[X_s X_t\right]=\prod_{i=1}^N
K_{1,H_i}\left[|s_i|^{2H_i}+|t_i|^{2H_i}-|t_i-s_i|^{2H_i}\right]
\end{eqnarray*}

\fin

This process also has an harmonizable representation, using the 
Fourier transform of the white noise in $\textbf{R}^N$ as in the previous 
paragraph. 

\begin{proposition} \label{fbsh}
For all $t=\left(t_i\right)$, consider the function $\phi_t$ such that for all 
$\xi=\left(\xi_i\right)$,
\[
\phi_t(u)=\prod_{m=1}^N \frac{e^{it_m \xi_m}-1}{|\xi_m|^{H_m+\frac{1}{2}}}
\]
The process defined by
\[
X_t=\hat{W}\left(\phi_t\right)=\int_{\textbf{R}^N} \prod_{m=1}^N \frac{e^{it_m \xi_m}-1}{|\xi_m|^{H_m+\frac{1}{2}}}\hat{\mathbb{W}}(d\xi)
\]
is a fractional Brownian sheet, up to a multiplicative constant.

\end{proposition}

\proof
As in the previous proposition, let's compute the covariance function of this
process.
\begin{eqnarray*}
E\left[X_s X_t\right]&=&\int_{\textbf{R}^N} \prod_{m=1}^N \frac{\left(e^{is_m
\xi_m}-1\right)\left(e^{-it_m \xi_m}-1\right)}{|\xi_m|^{2H_m+1}}.d\xi\\
&=&\prod_{m=1}^N \int_{\textbf{R}}\frac{\left(e^{is_m
\xi_m}-1\right)\left(e^{-it_m \xi_m}-1\right)}{|\xi_m|^{2H_m+1}}.d\xi_m\\
&=&\prod_{m=1}^N C_{1,H_m}\left[|s_m|^{2H_m}+|t_m|^{2H_m}-|t_m-s_m|^{2H_m}\right]
\end{eqnarray*}
using the same argument of the previous proposition.
\fin

\begin{remark}
The processes defined in propositions \ref{fbsmm} and \ref{fbsh} are proved to
have the same law. In fact, as a particular case of proposition \ref{mbsh}, they
are indistinguishable.
\end{remark}

\subsection{Stationarity of increments and self similarity}\label{increments}
Let us start by recalling the notion of increments in
$\textbf{R}_{+}^N$.\\ 
For a function $f : [0,1]^N\rightarrow\textbf{R}$ and $h\in\textbf{R}$,
one usually define the progressive difference in direction $\epsilon_i$ by
\begin{eqnarray*}
\Delta_{h,i} f(x)=\left\{\begin{array}{ll}f(x+h\epsilon_i)-f(x)&\textrm{ if
}x,x+h\epsilon_i\in [0,1]^N \\ 0&\textrm{ either}
\end{array}\right.
\end{eqnarray*}
and for $h\in\textbf{R}^N$ and $A=(i_1,\dots,i_k)$,
\[
\Delta_{h,A}f=\Delta_{h_{i_1},i_1}f \circ\dots\circ \Delta_{h_{i_k},i_k}f
\]
Despite the temptation to define the increments by
$X_t-X_s$ as in one dimension, it is better to set 
\begin{eqnarray} \label{defdelta}
\Delta X_{s,t}&=&\Delta_{t-s,(1,\dots,N)}X_s\nonumber\\
&=&\sum_{r\in\left\{0,1\right\}^N} (-1)^{N-\sum_l r_l}
X_{\left[s_i+r_i(t_i-s_i)\right]_i}
\end{eqnarray}
If there exists $i\in\left\{1,\dots,N\right\}$ such that $s_i=t_i$, 
we have $\Delta X_{s,t}=0$. Then, we consider
\[
I=\left\{i=1,\dots,N;\; s_i\neq t_i\right\}
\]
and
\[
\Delta_{t-s,I}X_s=\sum_{r\in\left\{0,1\right\}^{\#I}} (-1)^{\#I-\sum_l r_l}
X_{\left[s_i+r_i(t_i-s_i)\right]_{i\in I}}
\]

\subsubsection{Isotropic case}
In the isotropic case, the following extension of fBm's properties are
well known (see \cite{lindstr}).

\begin{proposition} \label{propisofbm}
Let $X=\left\{X_t; t\in\textbf{R}_{+}^N\right\}$ be a Levy fractional Brownian
motion. We have the two following properties for all $h\in\textbf{R}_{+}^N$
and $a>0$
\begin{eqnarray*}
X_{t+h}-X_h &\stackrel{(d)}{=}& X_t-X_0\\
X_{at} &\stackrel{(d)}{=}& a^H X_t
\end{eqnarray*}
where $\stackrel{(d)}{=}$ means equality of finite dimensional distributions.
\end{proposition}

\proof
For all $s$ and $t$ in $\textbf{R}_{+}^N$, we have
\begin{eqnarray*}
E\left[(X_{s+h}-X_h)(X_{t+h}-X_h)\right]&=&\frac{1}{2}\left(E\left[X_{s+h}-X_h\right]^2
+E\left[X_{t+h}-X_h\right]^2-E\left[X_{t+h}-X_{s+h}\right]^2\right)\\
&=&\frac{1}{2}\left(\|s\|^{2H}+\|t\|^{2H}-\|t-s\|^{2H}\right)\\
&=&E\left[X_s X_t\right]
\end{eqnarray*}

For self-similarity, we compute
\begin{eqnarray*}
E\left[X_{as}X_{at}\right]&=&\frac{1}{2}\left(E\left[X_{as}\right]^2+E\left[X_{at}\right]^2
-E\left[X_{at}-X_{as}\right]^2\right)\\
&=&\frac{1}{2}\left(\|as\|^{2H}+\|at\|^{2H}-\|at-as\|^{2H}\right)\\
&=&a^{2H} E\left[X_s X_t\right]\\
&=&E\left[a^H X_s\; a^H X_t\right]
\end{eqnarray*}

\fin
Proposition \ref{propisofbm} implies the stationarity of increments
(\ref{defdelta}).

\begin{proposition}\label{deltaiso}
The increments of Levy fractional Brownian are stationary, ie for all
$h\in\textbf{R}_{+}^N$
\[
\Delta X_{h,t+h} \stackrel{(d)}{=} \Delta X_{0,t}
\]
\end{proposition}

\proof
We fix $h\in\textbf{R}_{+}^N$ and write
\begin{eqnarray*}
\Delta X_{h,t+h}=\sum_{r\in\left\{0,1\right\}^N-\left\{0\right\}} (-1)^{N-\sum_l r_l}
\left(X_{\left[h_i+r_i t_i\right]_i}-X_h\right)
\end{eqnarray*}
then in the development of $E\left[\Delta X_{h,s+h}\Delta X_{h,t+h}\right]$,
we only have terms of the form
\begin{eqnarray*}
E\left[\left(X_{\left[h_i+r_i s_i\right]_i}-X_h\right)
\left(X_{\left[h_i+\rho_i t_i\right]_i}-X_h\right)\right]
=E\left[X_{\left[r_i s_i\right]_i}
X_{\left[\rho_i t_i\right]_i}\right]
\end{eqnarray*}
using the previous proposition.
Therefore we have
\[
E\left[\Delta X_{h,s+h}\Delta X_{h,t+h}\right]
=E\left[\Delta X_{0,s}\Delta X_{0,t}\right]
\]
\fin

\subsubsection{Non-isotropic case}
In the non-isotropic case, the properties of self-similarity and stationarity
of increments have been stated by L\'eger/Pontier (cf \cite{legerpontier}).
Here, we give another proof based on the covariance function rather than the 
moving average representation.

\begin{proposition}
Let $X=\left\{X_t; t\in\textbf{R}_{+}^N\right\}$ be a fractional Brownian sheet.
We have the two following properties for all $h\in\textbf{R}_{+}^N$
and $a>0$
\[
\Delta X_{h,t+h} \stackrel{(d)}{=} \Delta X_{0,t}
\]
\[
X_{at} \stackrel{(d)}{=} a^{\sum_i H_i}\; X_t
\]
\end{proposition}

\proof
We consider $N$ independent fBm $X^{(1)},\dots X^{(N)}$ of Hurst parameter 
$H_i$, and the process $Y=\left\{Y_t; t\in\textbf{R}_{+}^N\right\}$
such that $Y_t=\prod_{i=1}^N X^{(i)}_{t_i}$.
We can see easily that $X$ and $Y$ have the same covariance function. 
The same result follows for the increments 
$\left\{\Delta X_{h,t+h};t\in\textbf{R}_{+}^N\right\}$ and 
$\left\{\Delta Y_{h,t+h};t\in\textbf{R}_{+}^N\right\}$. 
As a consequence, from
\begin{eqnarray*}
\Delta Y_{h,t+h}&=&\sum_{r\in\left\{0,1\right\}^N} (-1)^{N-\sum_l r_l}
\prod_{i=1}^N X^{(i)}_{h_i+r_i t_i}\\
&=&\prod_{i=1}^N \left[X^{(i)}_{t_i+h_i}-X^{(i)}_{h_i}\right]
\end{eqnarray*}
then we have
\begin{eqnarray*}
E\left[\Delta X_{h,s+h} \Delta X_{h,t+h}\right]&=&E\left[\Delta Y_{h,s+h} 
\Delta Y_{h,t+h}\right]\\
&=&\prod_{i=1}^N \underbrace{E\left[\left(X^{(i)}_{s_i+h_i}-X^{(i)}_{h_i}\right)
\left(X^{(i)}_{t_i+h_i}-X^{(i)}_{h_i}\right)\right]}_{E\left[X^{(i)}_{s_i}X^{(i)}_{t_i}\right]}\\
&=&E\left[\Delta X_{0,s} \Delta X_{0,t}\right]
\end{eqnarray*}

For self-similarity, we compute for all $a>0$
\begin{eqnarray*}
E\left[X_{as}X_{at}\right]&=&\prod_{i=1}^{N}\frac{1}{2}\left((a s_i)^{2H_i}
+(a t_i)^{2H_i}-|a t_i-a s_i|^{2H_i}\right)\\
&=&a^{2\sum_i H_i}\prod_{i=1}^{N}\frac{1}{2}\left(s_i^{2H_i}
+t_i^{2H_i}-|t_i-s_i|^{2H_i}\right)\\
&=&a^{2\sum_i H_i}E\left[X_{s}X_{t}\right]\\
&=&E\left[a^{\sum_i H_i}\;X_{s}\; a^{\sum_i H_i}\;X_{t}\right]
\end{eqnarray*}

\fin

Therefore, we can conclude that both extensions of fBm satisfy the properties of 
self-similarity and stationarity of increments.

\section{The multifractional Brownian motion's case}
Once again, we can consider two different kinds of multi-parameter extension of 
mBm : isotropic and anisotropic extension.
Note, first of all, that mBm already has a multi-parameter extension. Indeed, 
the formulation of Benassi/Jaffard/Roux in \cite{benassi} was done for 
$t\in\textbf{R}^N$. We will see that it can be considered as an isotropic 
extension.

\subsection{Isotropic extension}
To define an isotropic extension of the mBm, the natural way is to substitute
the constant $H$ of the moving average representation of the Levy fractional
Brownian motion, with a function.

\begin{definition}
Let $H : \textrm{R}^N \rightarrow (0,1)$ be a measurable function.
The process $\left\{X_t; t\in\textbf{R}^N_{+}\right\}$ such that 
\begin{equation}\label{isomoy}
X_t=\int_{\textbf{R}^N} \left[\|t-u\|^{H(t)-\frac{N}{2}}-\|u\|^{H(t)-\frac{N}{2}}\right]\mathbb{W}(du)
\end{equation}
is called multifractional Brownian field.
\end{definition}

We will show that this process {\it is the same} as the process defined by 
Benassi/Jaffard/Roux. This result generalizes on the equivalence stated in the case 
$N=1$ in \cite{cohencov}.

\begin{proposition} \label{mbfieldmoy}
Let $H : \textrm{R}^N \rightarrow (0,1)$ be a measurable function.
The process defined by
\begin{equation}
X_t=\int_{\textrm{R}^N} \frac{e^{i\left<t,\xi\right>}-1}{\|\xi\|^{H(t)+\frac{N}{2}}}.\hat{\mathbb{W}}(d\xi)
\end{equation}
is indistinguishable, up to a multiplicative deterministic function,
from the process defined by (\ref{isomoy}).
This formulation is the harmonizable representation of the multifractional 
Brownian field.
\end{proposition}

\proof
First of all, let us compute the Fourier transform of the function
$\|.\|^{\alpha}$.
\begin{eqnarray*}
\left<\mathcal{T}\|.\|^{\alpha},\varphi\right>
&=&\left<\|.\|^{\alpha},\hat{\varphi}\right>\\
&=&\int_{\textbf{R}^N}\|t\|^{\alpha}\left(\int_{\textbf{R}^N}e^{-i<w,t>}\varphi(w).dw\right).dt
\end{eqnarray*}
we consider the change of variables
\begin{eqnarray*}
\textbf{R}^N\times\textbf{R}^N&\rightarrow&\textbf{R}^N\times\textbf{R}^N\\
(w,t)&\mapsto&(w,\lambda=\phi(t))
\end{eqnarray*}
where $\phi$ is the linear application which maps the canonic basis of $\textbf{R}^N$
to the orthonormal basis $\left(e_1=\frac{w}{\|w\|},e_2,\dots,e_N\right)$.
We get
\begin{eqnarray*}
\left<\mathcal{T}\|.\|^{\alpha},\varphi\right>
&=&\int_{\textbf{R}^N}\int_{\textbf{R}^N}\|\lambda\|^{\alpha}e^{i\lambda_1 \|w\|}\varphi(w).dw.d\lambda\\
&=&\int_{\textbf{R}^N}\int_{\textbf{R}^N}\frac{\|u\|^{\alpha}}{\|w\|^{\alpha}}
e^{-i u_1}\varphi(w)\frac{dw.du}{\|w\|^N}
\end{eqnarray*}
using the change of variables $(w,\lambda)\mapsto(w,u=\|w\|\lambda)$.
Then we have
\begin{eqnarray*}
\left<\mathcal{T}\|.\|^{\alpha},\varphi\right>
=\underbrace{\left(\int_{\textbf{R}^N}\|u\|^{\alpha}e^{-i u_1}.du\right)}_{\lambda_{\alpha}}
\int_{\textbf{R}^N}\frac{1}{\|w\|^{\alpha+N}}\; \varphi(w).dw
\end{eqnarray*}
Thus,  
\[
\mathcal{T}\|.\|^{\alpha}(w)=\frac{\lambda_{\alpha}}{\|w\|^{\alpha+N}}
\]

We use this result to calculate the Fourier transform of $\|t-.\|^{\alpha}-\|.\|^{\alpha}$.
We will use the following property : if $g(u)=f(u-\alpha)$ then $\hat{g}=e^{-i<\alpha,v>}\hat{f}(v)$.
\[
\mathcal{T}\left[\|t-.\|^{\alpha}-\|.\|^{\alpha}\right](v)=\left[e^{-i<t,v>}-1\right]\frac{\lambda_{\alpha}}{\|v\|^{\alpha+N}}
\]

We deduce from this
\[
\mathcal{T}\left[\|t-.\|^{H(t)-\frac{N}{2}}-\|.\|^{H(t)-\frac{N}{2}}\right](v)
=\lambda_{H(t)}\overline{\left(\frac{e^{i<t,v>}-1}{\|v\|^{H(t)+\frac{N}{2}}}\right)}
\]
and $\forall t\in\textbf{R}^N$, we have almost surely
\[
\int_{\textbf{R}^N} \left[\|t-u\|^{H(t)-\frac{N}{2}}-\|u\|^{H(t)-\frac{N}{2}}\right]\mathbb{W}(du)
=\lambda_{H(t)}\int_{\textrm{R}^N} \frac{e^{i\left<t,\xi\right>}-1}{\|\xi\|^{H(t)+\frac{N}{2}}}.\hat{\mathbb{W}}(d\xi) 
\]
using the fact we saw previously that the second integral is almost surely real.
Therefore, by an argument of continuity, the result follows.
\fin

This process is obviously a centered Gaussian process. It is thus of interest
to study its covariance function. The following proposition is an extension of
the case $N=1$ stated in \cite{cohen}.

\begin{proposition}
Let $\left\{X_t;t\in\textbf{R}_{+}^N\right\}$ be a multifractional Brownian
field. There exists a deterministic function 
$D^f_N:\textbf{R}\rightarrow \textbf{R}$ such that
the covariance function of $X$ can be written
\begin{equation} \label{covmbm2d}
E\left[X_s X_t\right]=D^f_N(H(s)+H(t))\left[\|s\|^{H(s)+H(t)}+\|t\|^{H(s)+H(t)}
-\|t-s\|^{H(s)+H(t)}\right]
\end{equation}

\end{proposition}

\proof
The easiest way to show this result is to use the harmonizable representation.
By definition of $\hat{\mathbb{W}}$, we have
\[
E\left[X_s X_t\right]=\int_{\textbf{R}^N}\frac{\left(e^{i<s,\xi>}-1\right) 
\left(e^{-i<t,\xi>}-1\right)}{\|\xi\|^{H(s)+H(t)+N}}.d\xi
\]
This integral has already been calculated for a Levy fractional Brownian motion
with a parameter $H=\frac{H(s)+H(t)}{2}$.
Then we have
\[
E\left[X_s X_t\right]=\underbrace{\left(\int_{\textbf{R}^N}\frac{1-e^{iu_1}}{\|u\|^{H(s)+H(t)+N}}.du\right)}
_{D^f_N(H(s)+H(t))}
\left[\|s\|^{H(s)+H(t)}+\|t\|^{H(s)+H(t)}-\|t-s\|^{H(s)+H(t)}\right]
\]
with $D^f_N(x)=\int_{\textbf{R}^N}\frac{1-e^{iu_1}}{\|u\|^{x+N}}.du$
\fin

\subsection{Non isotropic extension}
Another way to extend the multifractional Brownian motion for a set of index
included in $\textbf{R}^N_{+}$, is to copy the definition of the Brownian sheet.

\begin{definition}
Let $H:\textbf{R}^N_{+}\rightarrow (0,1)^N$ be a measurable function.
The process $\left\{X_t; t\in\textbf{R}_{+}^N\right\}$ such that
\[
X_t=\int_{\textbf{R}^N} \prod_{i=1}^N \left[|t_i-u_i|^{H_i(t)-\frac{1}{2}}-|u_i|^{H_i(t)-\frac{1}{2}}\right]\mathbb{W}(du)
\]
where $\mathbb{W}$ is the white noise, 
is called  multifractional Brownian sheet (mBs).
\end{definition}

As in the case of the isotropic extension, there also exists a harmonizable representation 
of the mBs.

\begin{proposition} \label{mbsh}
Let $H:\textbf{R}^N_{+}\rightarrow (0,1)^N$ be a measurable function.
For all $t=\left(t_i\right)_{i\in\left\{1;\dots;N\right\}}$, we consider 
the function $\phi_t$ such that for all $\xi=\left(\xi_i\right)$,
\[
\phi_t(u)=\prod_{m=1}^N \frac{e^{it_m \xi_m}-1}{|\xi_m|^{H_m(t)+\frac{1}{2}}}
\]
The process defined by
\[
X_t=\hat{W}\left(\phi_t\right)=\int_{\textbf{R}^N} \prod_{m=1}^N \frac{e^{it_m \xi_m}-1}{|\xi_m|^{H_m(t)+\frac{1}{2}}}\hat{\mathbb{W}}(d\xi)
\]
is indistinguishable, up to a multiplicative deterministic function,
from the process defined previously.
This formulation is the harmonizable representation of the multifractional 
Brownian sheet.
\end{proposition}

\proof
We have already seen that for each $m\in\left\{1,\dots,N\right\}$
\[
\mathcal{T}\left[|t_m-.|^{H_m(t)-\frac{1}{2}}-|.|^{H_m(t)-\frac{1}{2}}\right](\xi_m)
=\lambda_{H_m(t)}
\overline{\left(\frac{e^{it_m \xi_m}-1}{|\xi_m|^{H_m(t)+\frac{1}{2}}}\right)}
\]
Moreover, we compute
\begin{eqnarray*}
\mathcal{T}\left(\prod_{m=1}^N
\left[|t_m-.|^{H_m(t)-\frac{1}{2}}-|.|^{H_m(t)-\frac{1}{2}}\right]\right)&(\xi)&\\
=&&\int_{\textbf{R}^N} e^{-i<\xi,x>}\prod_{m=1}^N
\left[|t_m-x_m|^{H_m(t)-\frac{1}{2}}-|x_m|^{H_m(t)-\frac{1}{2}}\right].dx\\
=&&\int_{\textbf{R}^N} \prod_{m=1}^N e^{-i\xi_m x_m}\left[|t_m-x_m|^{H_m(t)-\frac{1}{2}}-|x_m|^{H_m(t)-\frac{1}{2}}\right].dx\\
=&&\prod_{m=1}^N
\mathcal{T}\left[|t_m-.|^{H_m(t)-\frac{1}{2}}-|.|^{H_m(t)-\frac{1}{2}}\right](\xi_m)
\end{eqnarray*}
Therefore
\[
\underbrace{\left(\prod_{i=1}^N \lambda_m(t)\right)}_{\lambda(t)}
\hat{W}\overline{\left(\prod_{m=1}^N \frac{e^{it_m.}-1}{|.|^{H_m(t)+\frac{1}{2}}}\right)}
=W\left(\prod_{m=1}^N \left[|t_m-.|^{H_m(t)-\frac{1}{2}}-|.|^{H_m(t)-\frac{1}{2}}\right]\right)
\]
We use the same arguments as in proposition \ref{mbfieldmoy} to conclude.
\fin

The following proposition shows that the covariance structure of multifractional
Brownian sheet, is a generalization of the fBs's one.

\begin{proposition}
Let $\left\{X_t;t\in\textbf{R}_{+}^N\right\}$ be a multifractional Brownian sheet. 
There exists a deterministic function $D^s:\textbf{R}^N\rightarrow\textbf{R}$ 
such that
\begin{equation} \label{covisombm}
E\left[X_s X_t\right]=D^s\left(H(s)+H(t)\right) \prod_{m=1}^N\left[|s_m|^{H_m(s)+H_m(t)}
+|t_m|^{H_m(s)+H_m(t)}-|t_m-s_m|^{H_m(s)+H_m(t)}\right]
\end{equation}
\end{proposition}

\proof
As usually, we use the harmonizable representation of the process
\begin{eqnarray*}
E\left[X_s X_t\right]&=&\int_{\textbf{R}^N} \prod_{m=1}^N \frac{\left(e^{is_m
\xi_m}-1\right)\left(e^{-it_m \xi_m}-1\right)}{|\xi_m|^{H_m(s)+H_m(t)+1}}.d\xi\\
&=&\prod_{m=1}^N \int_{\textbf{R}}\frac{\left(e^{is_m
\xi_m}-1\right)\left(e^{-it_m \xi_m}-1\right)}{|\xi_m|^{H_m(s)+H_m(t)+1}}.d\xi_m
\end{eqnarray*}
We remark that the factor corresponding to each $m$, is the covariance of a multifractional
Brownian motion, with has already been calculated.
Therefore we have
\begin{eqnarray*}
E\left[X_s X_t\right]=\prod_{m=1}^N D^f_1\left(H_m(s)+H_m(t)\right) \left[|s_m|^{H_m(s)+H_m(t)}
+|t_m|^{H_m(s)+H_m(t)}-|t_m-s_m|^{H_m(s)+H_m(t)}\right]
\end{eqnarray*}

\fin

\begin{remark} \label{prodtens}
The form of the previous covariance function gives the idea to consider
the process $Y=\left\{Y_t;t\in\textbf{R}_{+}^N\right\}$ defined from 
$N$ independent multifractional Brownian motions $X^{(i)}$ with parameter $H_i$ 
by
\[
Y_t=X^{(1)}_{t^{(1)}}\dots X^{(N)}_{t^{(N)}}
\]
Although $Y$ is not a Gaussian process, it is easily seen that it has the same
covariance function as a multifractional Brownian sheet. This remark will be
often used in the following.
\end{remark}

\section{Regularity}
A lot of properties are known about the regularity of the trajectories of Brownian
motion and fractional Brownian motion.
As we will see, in the case of the multi-parameter extension of the mBm, we have
to make some assumptions about the regularity of $H$ before studying the 
continuity of trajectories.
In the definitions of mBm (cf \cite{ayach} and \cite{benassi}), the function $H$ 
is  supposed to be H\"older continuous.

\subsection{Continuity of the two extensions}
We first recall the Kolmogorov's criterion.
\begin{theorem}[Kolmogorov]
Consider a process $X=\left\{X_t;t\in\textbf{R}_{+}^N\right\}$ such that
there exists $C>0$, $p>0$ and $\gamma>N$
\begin{equation}
\forall s,t\in\textbf{R}_{+}^N,\; E\left[|X_t-X_s|^p\right]\leq C \|t-s\|^{\gamma}
\end{equation}
Then, there exists a modification $Y=\left\{Y_t;t\in\textbf{R}_{+}^N\right\}$
of $X$ that is H\"older continuous of any order $q\in (0,\frac{\gamma-N}{p})$.
\end{theorem}

As usually, the quantity $E\left[|X_t-X_s|^2\right]$ is studied for
$s,t\in[a,b]$ where $a\preceq b$ and then, a patching argument is used to
extend to $s,t\in\textbf{R}_{+}^N$.

\subsubsection{Isotropic case}

\begin{lemma}\label{lemD}
For all $\eta$ and $\mu$ such that $0<\eta<\mu<1$,
the multiplicative factor $D^f_N$ of covariance function in (\ref{covisombm}), 
is positive and belongs to $C^{\infty}\left([\eta,\mu]\right)$.
Moreover, its order $n$ derivative is given by
\begin{equation}\label{derD}
{D^f_N}^{(n)}(x)=\int_{\textbf{R}^N}\frac{1-e^{iu_1}}{\|u\|^{x+N}}\ln^n\frac{1}{\|u\|}.du 
\end{equation}
\end{lemma}

\proof
As the integral of a positive function, $D^f_N$ is positive. By an argument of uniform
convergence of integrals (\ref{derD}) on $[\eta,\mu]$, $D^f_N$ is 
$C^{\infty}\left([\eta,\mu]\right)$ and the derivatives are obtained by
derivations of the integrand.
\fin

\begin{proposition} \label{propcov}
For all $s,t\in [a,b]$, we have
\begin{eqnarray} \label{covisoasymp}
\frac{1}{2}E\left[X_t-X_s\right]^2&=&D\left[H(s)+H(t)\right]\times\|t-s\|^{H(s)+H(t)}\nonumber\\
&&+\frac{1}{2}\left[\frac{\partial^2\varphi}{\partial
x^2}\left(H(s)+H(t);\|s\|\right)+\frac{\partial^2\varphi}{\partial
x^2}\left(H(s)+H(t);\|t\|\right)\right]\times\left(H(t)-H(s)\right)^2\nonumber\\
&&+O_{a,b}\left[\left(H(t)-H(s)\right)\left(\|t\|-\|s\|\right)\right]
+o_{a,b}\left(H(t)-H(s)\right)^2
\end{eqnarray}
where $\varphi(x,y)=D(x)y^x$.
\end{proposition}

\proof
Using the covariance function of the multifractional Brownian field,
we have
\begin{eqnarray} \label{kolmcoviso}
\frac{1}{2}E\left[|X_s-X_t|^2\right]&=&D\left[2H(s)\right]\; \|s\|^{2H(s)}-D\left[H(s)+H(t)\right]\; \|s\|^{H(s)+H(t)}\nonumber\\
&&+D\left[2H(t)\right]\; \|t\|^{2H(t)}-D\left[H(s)+H(t)\right]\; \|t\|^{H(s)+H(t)}\nonumber\\
&&+D\left[H(s)+H(t)\right]\; \|t-s\|^{H(s)+H(t)}
\end{eqnarray}
We have to get a second order expansion of this expression.\\
We introduce the function $\varphi$ defined by
\[
\varphi(x,y)=D(x)y^x
\]
We can write
\begin{eqnarray} \label{covisophi}
\frac{1}{2}E\left[|X_s-X_t|^2\right]&=&\varphi(2H(s),\|s\|)-\varphi(H(s)+H(t),\|s\|)\nonumber\\
&&+\varphi(2H(t),\|t\|)-\varphi(H(s)+H(t),\|t\|)\nonumber\\
&&+D\left[H(s)+H(t)\right]\; \|t-s\|^{H(s)+H(t)}
\end{eqnarray}
We use the second order expansion
\begin{eqnarray*}
\varphi(2H(s),\|s\|)-\varphi(H(s)+H(t),\|s\|)
&=&\left(H(s)-H(t)\right)\times\frac{\partial\varphi}{\partial x}\left(H(s)+H(t),\|s\|\right)\\
&&+\frac{\left(H(s)-H(t)\right)^2}{2}\times\frac{\partial^2\varphi}{\partial x^2}\left(H(s)+H(t),\|s\|\right)\\
&&+o_{a,b}\left(H(s)-H(t)\right)^2
\end{eqnarray*}
An inversion of roles between $s$ and $t$ provides the expansion of \\
\[
\varphi(2H(t),\|t\|)-\varphi(H(s)+H(t),\|t\|)
\]
Then (\ref{covisophi}) becomes
\begin{eqnarray*}
\frac{1}{2}E\left[|X_s-X_t|^2\right]&=&\left(H(t)-H(s)\right)\times
\left[\frac{\partial\varphi}{\partial x}\left(H(s)+H(t),\|t\|\right)
-\frac{\partial\varphi}{\partial x}\left(H(s)+H(t),\|s\|\right)\right]\\
&&+\frac{\left(H(t)-H(s)\right)^2}{2}\times\left[
\frac{\partial^2\varphi}{\partial x^2}\left(H(s)+H(t),\|s\|\right)
+\frac{\partial^2\varphi}{\partial x^2}\left(H(s)+H(t),\|t\|\right)\right]\\
&&+D\left[H(s)+H(t)\right]\; \|t-s\|^{H(s)+H(t)}
+o_{a,b}\left(H(t)-H(s)\right)^2
\end{eqnarray*}
Since 
\[\left(H(t)-H(s)\right)\times
\left[\frac{\partial\varphi}{\partial x}\left(H(s)+H(t),\|t\|\right)
-\frac{\partial\varphi}{\partial x}\left(H(s)+H(t),\|s\|\right)\right]
\]
is
$O_{a,b}\left[\left(H(t)-H(s)\right)\left(\|t\|-\|s\|\right)\right]$, the result
follows.
\fin

\begin{corollary}\label{corisoasymp}
For all $s,t\in [a,b]$, we have
\begin{eqnarray} \label{covisoasympsym}
\frac{1}{2}E\left[X_t-X_s\right]^2&=&D\left[2H(t)\right]\times\|t-s\|^{2H(t)}\nonumber\\
&&+\frac{\partial^2\varphi}{\partial x^2}\left(2H(t);\|t\|\right)\times\left(H(t)-H(s)\right)^2\nonumber\\
&&+o_{a,b}\left(H(t)-H(s)\right)^2+o_{a,b}\left(\|t-s\|^{2H(t)}\right)
\end{eqnarray}
where $\varphi(x,y)=D(x)y^x$.
\end{corollary}

\proof
Using the expansion of $D\left[H(s)+H(t)\right]$ and
\begin{eqnarray*}
\|t-s\|^{H(s)+H(t)}=\|t-s\|^{2H(t)}-\left(H(t)-H(s)\right)\|t-s\|^{2H(t)}\ln\|t-s\|
+o_{a,b}\left(H(t)-H(s)\right)^2
\end{eqnarray*}
we get
\begin{eqnarray} \label{cor1}
D\left[H(s)+H(t)\right]\times\|t-s\|^{H(s)+H(t)}&=&D\left[2H(t)\right]\times
\|t-s\|^{2H(t)}\\
&&+o_{a,b}\left(\|t-s\|^{2H(t)}\right)+o_{a,b}\left(H(t)-H(s)\right)^2\nonumber
\end{eqnarray}
Moreover as $H(t)<1$ for all $t\in [a,b]$, we have $\epsilon=1-H(t)>0$ and
\begin{eqnarray*}
2\left(H(t)-H(s)\right)\left(\|t\|-\|s\|\right)&=&
2\left(H(t)-H(s)\right)\left(\|t\|-\|s\|\right)^{\frac{\epsilon}{2}}\times
\left(\|t\|-\|s\|\right)^{1-\frac{\epsilon}{2}}\\
&\leq&\left(H(t)-H(s)\right)^2 \left(\|t\|-\|s\|\right)^{\epsilon}
+\left(\|t\|-\|s\|\right)^{2-\epsilon}
\end{eqnarray*}
that implies
\begin{eqnarray} \label{cor2}
\left(H(t)-H(s)\right)\left(\|t\|-\|s\|\right)=
o_{a,b}\left(H(t)-H(s)\right)^2+o_{a,b}\left(\|t-s\|^{2H(t)}\right)
\end{eqnarray}
We conclude by (\ref{covisoasymp}), (\ref{cor1}) and (\ref{cor2}) 
using first order expansion of 
$\frac{\partial^2\varphi}{\partial x^2}$ in $x$ and $y$.
\fin

Using the continuity of $D$, $D'$ and $D''$, we can state from the previous
proposition
\begin{corollary} \label{kolmiso1}
There exist positive constants $K$ and $L$ such that
\begin{equation}\label{majcoviso}
\forall s,t\in [a,b];\; E\left[X_t-X_s\right]^2\leq K\;\|t-s\|^{2H(t)} + L\;\left|H(t)-H(s)\right|^2
\end{equation}

\end{corollary}

\begin{corollary} \label{kolmiso}
Suppose $H$ is $\beta$-H\"older continuous.
There exists a constant $M$ such that
\begin{equation}
\forall s,t\in [a,b];\; E\left[X_t-X_s\right]^2\leq M\; \|t-s\|^{2(\beta\wedge H(t))}
\end{equation}

\end{corollary}

\subsubsection{Non-isotropic case} \label{parkolmdrap}

\begin{lemma} \label{kolmdrap1}
There exists positive constants $K$ and $L$ such that
\begin{equation}
\forall s,t\in [a,b];\; E\left[|X_t-X_s|^2\right]\leq K\;
\|t-s\|^{2\min_i H_i(t)}
+L\;\|H(t)-H(s)\|^2
\end{equation}
\end{lemma}

\proof
By remark \ref{prodtens},
we have
\begin{eqnarray*}
E\left[X_s-X_t\right]^2&=&E\left[\prod_{i=1}^N X^{(i)}_{s^{(i)}} - \prod_{i=1}^N X^{(i)}_{t^{(i)}}\right]^2\\
&=&E\left[\left(\prod_{i=1}^N X^{(i)}_{s^{(i)}} - X^{(1)}_{t^{(1)}}\prod_{i>1}X^{(i)}_{s^{(i)}}\right)
+\left(X^{(1)}_{t^{(1)}}\prod_{i>1}X^{(i)}_{s^{(i)}} -
X^{(1)}_{t^{(1)}}X^{(2)}_{t^{(2)}}\prod_{i>2}X^{(i)}_{s^{(i)}}\right)\right.\\
&&+\dots+
\left.\left(\left(\prod_{i=1}^{N-1} X^{(i)}_{t^{(i)}}\right)X^{(N)}_{s^{(N)}} -
\prod_{i=1}^N X^{(i)}_{t^{(i)}}\right)\right]^2
\end{eqnarray*}
By the inequality of convexity $\left(\sum a_i\right)^2\leq n\sum a_i^2$,
we get
\begin{eqnarray*}
E\left[X_s-X_t\right]^2\leq N\left\{E\left[\prod_{i>1}X^{(i)}_{s^{(i)}}\right]^2 E\left[X^{(1)}_{s^{(1)}}-X^{(1)}_{t^{(1)}}\right]^2
+\cdots+E\left[X^{(N)}_{s^{(N)}}-X^{(N)}_{t^{(N)}}\right]^2 E\left[\prod_{i=1}^{N-1}X^{(i)}_{s^{(i)}}\right]^2\right\}
\end{eqnarray*}
Since there exists a constant $M=M_{a,b}$ such that
\[
\forall t\in [a,b],\forall i;\; E\left[X^{(i)}_{t^{(i)}}\right]^2\leq M
\]
we get
\begin{eqnarray}\label{majinterm}
E\left[X_s-X_t\right]^2\leq N M^{n-1} \sum_{i=1}^N E\left[X^{(i)}_{s^{(i)}}-X^{(i)}_{t^{(i)}}\right]^2
\end{eqnarray}
Using
\[
E\left[X^{(i)}_{s^{(i)}}-X^{(i)}_{t^{(i)}}\right]^2\leq 
K_i |s^{(i)}-t^{(i)}|^{2H_i(t)}+L_i\left(H_i(s)-H_i(t)\right)^2
;\;\forall i=1,\dots,N
\]
(\ref{majinterm}) implies
\begin{eqnarray*}
E\left[X_t-X_s\right]^2\leq N M^{n-1}\left[\left(\sum_{i=1}^N K_i\right) 
\|t-s\|^{2\min_i H_i(t)}+\left(\sum_{i=1}^N L_i\right)
\|H(t)-H(s)\|^2 \right]
\end{eqnarray*}
\fin

\begin{corollary}
Suppose $H$ is $\beta$-H\"older continuous.
There exists a positive constant $M$ such that
\begin{equation}
\forall s,t\in [a,b];\; E\left[X_t-X_s\right]^2\leq M
\|t-s\|^{2\left(\beta\wedge\min_i H_i(t)\right)}
\end{equation}

\end{corollary}

\subsubsection{Existence of a continuous modification} \label{kolmconcl}
In both isotropic and anisotropic cases, under H\"older regularity assumptions
for $H$, we have an inequality
\[
E\left[X_t-X_s\right]^2\leq K\|t-s\|^{\alpha}
\]
But to use the Kolmogorov criterion, we need to have $\alpha>N$.\\
As the random variable $X_t-X_s$ is Gaussian, we can write, for each integer $n$
\[
E\left[X_t-X_s\right]^{2n}\leq \lambda_n K \|t-s\|^{n.\alpha}
\]
and choose $n$ such that $n.\alpha > N$.\\ 
We conclude by a classical patching argument.
For $a$ and $b$, Kolmogorov's theorem gives a continuous process
$Y^{a,b}=\left\{Y^{a,b}_t;t\in [a,b]\right\}$.
Consider $a'$ and $b'$ such that $[a,b]\subset [a',b']$. The processes
$Y^{a,b}$ and $Y^{a',b'}$ coincide on $[a,b]$. Thus
\[
\forall t\in [a,b];\; P\left\{Y^{a,b}_t=Y^{a',b'}_t\right\}=1
\]
and, by continuity
\[
P\left\{Y^{a,b}_t=Y^{a',b'}_t;\forall t\in [a,b]\right\}=1
\]
Then we can define a process $Y$ on $\textbf{R}_{+}^N$ who coincides with
$Y^{a,b}$ on $[a,b]$ and we can seen easily that this process is continuous.

\subsection{H\"older exponents}
The notion of H\"older function is well known. It is interesting to consider
a localized version of this notion.

For the paths of a process $X$, one usually define two kinds of exponent (see
\cite{ayach}, \cite{GmBm}):
\begin{itemize}
\item the pointwise H\"older exponent
\begin{eqnarray*}
\alpha(t_0)&=&\sup\left\{\alpha; \lim_{h\rightarrow 0}\frac{|X_{t_0+h}-X_{t_0}|}{\|h\|^{\alpha}}=0\right\}\\
&=&\sup\left\{\alpha;\; \limsup_{\rho\rightarrow 0}\frac{\sup_{s,t\in B(t_0,\rho)}|X_t-X_s|}{\rho^{\alpha}}<\infty\right\}
\end{eqnarray*}

\item the local H\"older exponent
\begin{eqnarray*}
\tilde{\alpha}(t_0)=\sup\left\{\alpha;\; \limsup_{\rho\rightarrow 0}\sup_{s,t\in B(t_0,\rho)}\frac{|X_t-X_s|}{\|t-s\|^{\alpha}}<\infty\right\}
\end{eqnarray*}
\end{itemize}

We can see easily that for all $t_0$, we have
\begin{equation}\label{locinfpoint}
\tilde{\alpha}(t_0)\leq \alpha(t_0)
\end{equation}
A study of these exponents, in the case of $1$D mBm, is made in \cite{GmBm}.
\begin{remark}
If $H$ is $\beta$-H\"older continuous, then the local
H\"older exponent $\tilde{\beta}(t)$ of $H$ at every point is not smaller than 
$\beta$.
 
Conversely, suppose that the local H\"older exponent of $H$ at every point
of a compact $[a,b]$ is positive. Then $H$ is $\beta$-H\"older continuous 
on $[a,b]$  with
$\beta=\inf_{t\in[a,b]}\tilde{\beta}(t)$.
\end{remark}

In the following, we suppose that $H$ admits positive local H\"older exponent 
$\tilde{\beta}(t_0)$ at every point $t$.

\begin{proposition} \label{propexpiso}
Let $X=\left\{X_t; t\in\textbf{R}_{+}^N\right\}$ be a multifractional Brownian field.
For all $t_0\in\textbf{R}_{+}^N$, the local H\"older exponent of $X$ at $t_0$ is
almost surely given by
\begin{eqnarray}
\tilde{\alpha}(t_0)=\tilde{\beta}(t_0)\wedge H(t_0)
\end{eqnarray}
and the pointwise H\"older exponent of $X$ at $t_0$ satisfies almost surely
\begin{eqnarray}
\alpha(t_0)=\beta(t_0)\wedge H(t_0)
\end{eqnarray}
where $\beta(t_0)$ and $\tilde{\beta}(t_0)$ denote the pointwise and local
H\"older exponents of $H$ at $t_0$.
\end{proposition}

As a consequence of this result, if $H$ satisfies 
\[
\forall t\in\textbf{R}_{+}^N;\;\beta(t)<H(t)
\]
the H\"older regularity of multifractional Brownian field of parameter function
$H$ is given by the regularity of $H$ (and not by the value of $H$). This point
is developed in \cite{ehjlv}.

\begin{proposition} \label{propexpdrap}
Let $X=\left\{X_t; t\in\textbf{R}_{+}^N\right\}$ be a multifractional Brownian
sheet.
For all $t_0\in\textbf{R}_{+}^N$, the local H\"older exponent of $X$ at $t_0$ is
almost surely given by
\begin{eqnarray}
\tilde{\alpha}(t_0)=\tilde{\beta}(t_0)\wedge\min_i H_i(t_0)
\end{eqnarray}
and the pointwise H\"older exponent of $X$ at $t_0$ satisfies almost surely
\begin{eqnarray}
\alpha(t_0)=\beta(t_0)\wedge\min_i H_i(t_0)
\end{eqnarray}
where $\beta(t_0)$ and $\tilde{\beta}(t_0)$ denote the pointwise and local
H\"older exponents of $H$ at $t_0$.
\end{proposition}
The proofs of propositions \ref{propexpiso} and \ref{propexpdrap} are detailed
in the three following paragraphs.

\subsubsection{Lower bound for the local H\"older exponent} \label{lblhe}
A lower bound for the local H\"older exponent is directly given by 
Kolmogorov's theorem.
Indeed, for $X$ a multifractional Brownian field 
or a multifractional Brownian sheet indexed by $[a,b]$, for all $n\in\textbf{N}$, 
there exists $\lambda_n>0$ such that
\[
E\left[X_t-X_s\right]^{2n}\leq \lambda_n  \|t-s\|^{n.\alpha}
\]
with $\alpha=2\inf_{[a,b]}(\tilde{\beta}\wedge H)$ or
$\alpha=2\inf_{[a,b]}(\tilde{\beta}\wedge\min_i H_i)$.\\
Kolmogorov's theorem states that there exists a modification of $X$,
which is $q$-H\"older continuous for all $q\in(0,\frac{\alpha}{2}-\frac{N}{2n})$.
Then, for all $t_0\in\textbf{R}_{+}^N$ and all $a,b\in\textbf{R}_{+}^N$
such that $a\prec b$ and $t_0\in(a,b)$, 
we have
\[
\forall n\in\textbf{N};\;
\tilde{\alpha}(t_0)\geq \frac{\alpha}{2}-\frac{N}{2n}
\]
and therefore, taking the limit $n\rightarrow\infty$
\[
\tilde{\alpha}(t_0)\geq \frac{\alpha}{2}
\]
As $H$ is continuous, we can take the limit $(a,b)\rightarrow(t_0,t_0)$
and we get
\begin{itemize}
\item in the isotropic case,
\begin{equation}
\tilde{\alpha}(t_0)\geq \tilde{\beta}(t_0)\wedge H(t_0)
\end{equation}
\item in the non-isotropic case,
\begin{equation}
\tilde{\alpha}(t_0)\geq \tilde{\beta}(t_0)\wedge\min_i H_i(t_0)
\end{equation}
\end{itemize}

\subsubsection{Lower bound for the pointwise H\"older exponent}
By (\ref{locinfpoint}), paragraph \ref{lblhe} provides a lower bound for the 
pointwise H\"older exponent. However, it can be improved in the case 
$\tilde{\beta}(t_0)<\beta(t_0)$.\\
Let $X=\left\{X_t;t\in\textbf{R}_{+}^N\right\}$ be a multifractional Brownian
field.
By corollary \ref{kolmiso1}, there exist positive constants $K$ and $L$ 
such that for all $s,t\in\textbf{R}_{+}^N$,
\[
E\left[X_t-X_s\right]^2\leq K\;\|t-s\|^{2H(t)}
+L\;|H(t)-H(s)|^2
\]
and by corollary \ref{kolmiso}, there exists positive
constants $\alpha$ and $M$ such that
\[
\forall s,t\in[a,b];\; E\left[X_t-X_s\right]^2\leq M \|t-s\|^{\alpha}
\]
Therefore, using Kolmogorov's criterion, there exists a modification of $X$, which
is $\nu$-H\"older continuous for all $\nu\in\left]0,\frac{\alpha}{2}\right[$. In
the following, we consider such a $\nu$ with $\frac{1}{\nu}\in\textbf{N}$\\
For all $\epsilon>0$, there exist $\rho_0>0$
and $M>0$ such that for all $\rho<\rho_0$ and all $t\in B(t_0,\rho)$
\[
E\left[\frac{X_t-X_{t_0}}{\rho^{\beta(t_0)\wedge
H(t_0)-\epsilon}}\right]^2
\leq M\;\rho^{\epsilon}
\]
Then, setting $\gamma=\beta(t_0)\wedge H(t_0)-\epsilon$,
for all $p\in\textbf{N}^{*}$
\begin{eqnarray*}
P\left\{|X_t-X_{t_0}|>\rho^{\gamma}\right\}
\leq E\left[\frac{X_t-X_{t_0}}{\rho^{\gamma}}\right]^{2p} \leq M_p\;\rho^{p\epsilon}
\end{eqnarray*}
Let $\rho=2^{-n}$ and for all $m\in\textbf{N}$,
\[
D_m=\left\{t_0+k.2^{-(n+m)};\;k\in\left\{0,\pm 1,\dots,\pm 2^{m}\right\}^N\right\}
\]
let us compute
\begin{eqnarray*}
&&P\left\{\max_{k\in\left\{\pm 1,\dots,\pm 2^m\right\}^N}
\frac{|X_{t_0+k.2^{-(m+n)}}-X_{t_0}|}{2^{-\gamma n}} > 1\right\}\\
&&\leq \sum_{k\in\left\{\pm 1,\dots,\pm 2^m\right\}^N}
P\left\{|X_{t_0+k.2^{-(n+m)}}-X_{t_0}|
>2^{-\gamma n}\right\}\\
&&\leq M_p\;2^{(m+1)N}\; 2^{-p\epsilon n}
\end{eqnarray*}
Let us take $m=\frac{1+\lfloor\gamma\rfloor}{\nu}n=\kappa n$ and $p\in\textbf{N}$ such that
$N\frac{1+\lfloor\gamma\rfloor}{\nu}-p\epsilon<0$.
By the Borel-Cantelli lemma, there
exists a finite random variable $n^{*}$ such that almost surely,
\begin{equation} \label{maxmaj}
\forall n\geq n^{*};\;
\max_{k\in\left\{0,\dots,\pm 2^{\kappa n}\right\}^N}
|X_{t_0+k.2^{-(1+\kappa)n}}-X_{t_0}|
\leq 2^{-\gamma n}
\end{equation}
From (\ref{maxmaj}), we show that, almost surely, for all $m\in\textbf{N}$, we
have
\begin{eqnarray}\label{majrec}
\forall t\in D_m;\;|X_t-X_{t_0}|\leq C\; 2^{-\gamma n}
\end{eqnarray}
\begin{itemize}
\item if $0\leq m\leq\kappa n$, (\ref{majrec}) follows directly from (\ref{maxmaj})
\item if $m>\kappa n$,
for $t\in D_{m}$, let 
\[
C^{\kappa n}_{t_0,t}=\left\{x\in D_{\kappa n};\forall i, 
(t_0)_i\leq x_i\leq t_i\right\}
\]
Then consider 
$\hat{t}\in B(t,2^{-(1+\kappa)n})\cap C^{\kappa n}_{t_0,t}$.\\
As the paths of $X$ are $\nu$-H\"older continuous, we have
\begin{eqnarray*}
|X_{\hat{t}}-X_t|\leq \tilde{C} 2^{-\nu(1+\kappa) n} \leq \tilde{C} 2^{-\gamma n}
\end{eqnarray*}
and by (\ref{maxmaj}),
\[
|X_{\hat{t}}-X_{t_0}|\leq 2^{-\gamma n}
\]
Using the triangular inequality, the result follows.
\end{itemize}
Therefore, (\ref{majrec}) leads to
\begin{eqnarray*}
\forall m\in\textbf{N}; \forall s,t\in D_m;\;
|X_t-X_s|\leq 2C\;2^{-\gamma n}
\end{eqnarray*}
Using the continuity of $X$ and $m\rightarrow +\infty$, we get
\[
\sup_{s,t\in B(t_0,2^{-n})}|X_t-X_s|\leq 2C\;2^{-\gamma n}
\]
and therefore, almost surely,
\begin{equation} \label{limfonct}
\limsup_{\rho\rightarrow 0}\sup_{s,t\in B(t_0,\rho)}
\frac{|X_t-X_s|}{\rho^{\gamma}}<+\infty
\end{equation}
By (\ref{limfonct}), for all $\epsilon>0$, almost surely
\[
\alpha(t_0)\geq \beta(t_0)\wedge H(t_0)-\epsilon
\]
Taking $\epsilon\in\textbf{Q}_{+}$, we have almost surely
\begin{equation}
\alpha(t_0)\geq \beta(t_0)\wedge H(t_0)
\end{equation}
For a multifractional Brownian sheet $X$, by lemma \ref{kolmdrap1}, we get in
the same way that, almost surely
\begin{equation}
\alpha(t_0)\geq \beta(t_0)\wedge H_i(t_0)
\end{equation}
for all $i=1,\dots,N$.

\subsubsection{Upper bound for the pointwise H\"older exponent}
The main result getting the upper bound for the H\"older exponents, is 
the following lemma, a direct
consequence of proposition \ref{propcov} using continuity of $D$, $D'$ and
$D''$.
\begin{lemma} \label{upisolem}
Let  $X=\left\{X_t; t\in\textbf{R}_{+}^N\right\}$ be a multifractional Brownian
field. For all $[a,b]\subset\textbf{R}_{+}^N$, there exist positive
constants $k_1$, $k_2$, $l_1$, $l_2$ such that
\begin{eqnarray}
\forall s,t\in[a,b];&&E\left[X_t-X_s\right]^2\geq
k_1\;\|t-s\|^{2H(t)}-l_1\left(H(t)-H(s)\right)^2 \label{miniso1}\\
&&E\left[X_t-X_s\right]^2\geq k_2\left(H(t)-H(s)\right)^2 \label{miniso2}
-l_2\;\|t-s\|^{2H(t)}
\end{eqnarray}
\end{lemma}

\proof
We only have to study the multiplicative factors of $\|t-s\|^{2H(t)}$ and
$\left(H(t)-H(s)\right)^2$ in (\ref{covisoasymp}) 
\begin{itemize}
\item Let 
$k_1=\inf_{t\in[a,b]}D\left[2H(t)\right]$
and
$l_2=\sup_{t\in[a,b]}D\left[2H(t)\right]$.\\
By continuity of $t\mapsto D\left[2H(t)\right]$ on the compact $[a,b]$
and as the function $D$ is positive (lemma \ref{lemD}), for all $t\in[a,b]$
\[
0<k_1\leq D\left[2H(t)\right]\leq l_2<+\infty
\]

\item and let
\[
\Phi(t)=\|t\|^{2H(t)}\times\left\{D\left[2H(t)\right]\ln^2\|t\|
-2D'\left[2H(t)\right]\ln\|t\|+D''\left[2H(t)\right]\right\}
\]
By lemma \ref{lemD},
\begin{eqnarray*}
\Phi(t)=\int_{\textbf{R}^N} \underbrace{\frac{1-e^{iu_1}}{\|u\|^{2H(t)+N}}
\left(\ln\|t\|-\ln\|u\|\right)^2}_{\geq 0}.du&&
\end{eqnarray*}
Let
$k_2=\inf_{t\in[a,b]}\Phi(t)$
and
$l_1=\sup_{t\in[a,b]}\Phi(t)$.\\
As previously, for all $t\in[a,b]$, we have
\[
0<k_2\leq\Phi(t)\leq l_1<+\infty
\]

\end{itemize}

\fin

\begin{lemma}
Let  $X=\left\{X_t; t\in\textbf{R}_{+}^N\right\}$ be a multifractional Brownian
sheet. 
For all $[a,b]\subset\textbf{R}_{+}^N$, there exist positive
constants $k_1$, $k_2$, $l_1$, $l_2$ such that
\begin{eqnarray}
\forall s,t\in[a,b];
&&t-s\in\textbf{R}_{+}.\epsilon_i\nonumber\\
&&E\left[X_t-X_s\right]^2\geq
k_1\;\|t-s\|^{2H_i(t)}-l_1\left(H_i(t)-H_i(s)\right)^2\label{mindrap1}\\
&&E\left[X_t-X_s\right]^2\geq
k_2\left(H_i(t)-H_i(s)\right)^2-l_2\;\|t-s\|^{2H_i(t)}\label{mindrap2}
\end{eqnarray}
\end{lemma}

\proof
For all $s$, $t$ such that 
$t-s\in\textbf{R}_{+}.\epsilon_i$, using lemma \ref{upisolem},
we have
\begin{eqnarray*}
E\left[X_t-X_s\right]^2&=&E\left[X^{(i)}_{t^{(i)}}-X^{(i)}_{s^{(i)}}\right]^2
\prod_{j\neq i} E\left[X^{(j)}_{t^{(j)}}\right]^2\\
&\geq&k_1\; |t_i-s_i|^{2H_i(t)}-l_1\left(H_i(t)-H_i(s)\right)^2
\end{eqnarray*}
and
\begin{eqnarray*}
E\left[X_t-X_s\right]^2\geq
k_2\left(H_i(t)-H_i(s)\right)^2-l_2\;|t_i-s_i|^{2H_i(t)}
\end{eqnarray*}

\fin

From this result, the upper bound for the pointwise exponent is a 
consequence of the following lemma whose proof is the same as the case $N=1$
(see \cite{ayach})

\begin{lemma} \label{majloclem}
Let $X=\left\{X_t;t\in\textbf{R}_{+}^N\right\}$ be a Gaussian process.
Assume there exists $\mu\in (0,1)$ such that for all $\epsilon>0$, there exist
a sequence $\left(h_n\right)_{n\in\textbf{N}}$ of
$\left(\textbf{R}_{+}^N\right)^{*}$ converging to 0, and a constant $c>0$
such that
\[
\forall n\in\textbf{N};\;
E\left[X_{t+h_n}-X_t\right]^2\geq c\|h_n\|^{2\mu+\epsilon}
\]
Then we have almost surely
\[
\alpha(t)\leq\mu
\]

\end{lemma} 

Let $X=\left\{X_t;t\in\textbf{R}^N_{+}\right\}$ be a multifractional Brownian field
(resp. multifractional Brownian sheet).
Let $\beta(t_0)$ be the pointwise H\"older exponent of $H$ at $t_0$.\\
We consider the two cases :
\begin{itemize}
\item if $H(t_0)<\beta(t_0)$ (resp. $H_i(t_0)<\beta(t_0)$), 
by definition of $\beta(t_0)$, we have
\begin{eqnarray*}
\lim_{h\rightarrow 0}\frac{\|H(t_0+h)-H(t_0)\|}{\|h\|^{H(t_0)}}=0
\end{eqnarray*}
Hence, by (\ref{miniso1}) (resp. (\ref{mindrap1})), there exists a positive constant $C$ such that
\begin{eqnarray*}
E\left[X_{t_0+h}-X_{t_0}\right]^2\geq C \|h\|^{2H(t_0)}
\end{eqnarray*}
Then, by lemma \ref{majloclem}
\begin{equation}\label{upponct1}
\alpha(t_0)\leq H(t_0)\textrm{ (resp. }H_i(t_0)\textrm{ )}
\end{equation}

\item if $H(t_0)>\beta(t_0)$ (resp.$H_i(t_0)>\beta(t_0)$),
we consider $\alpha\in (\beta(t_0);H(t_0))$ 
(resp. $\alpha\in (\beta(t_0);H_i(t_0))$).
There exists a positive constant $C$ and a sequence
$\left(h_n\right)_{n\in\textbf{N}}$ converging to $0$ such that
\begin{eqnarray*}
\forall n\in\textbf{N};\;\|H(t_0+h_n)-H(t_0)\|>C\|h_n\|^{\alpha}
\end{eqnarray*}
Then, by (\ref{miniso2}) (resp. (\ref{mindrap2}))
\begin{eqnarray*}
\forall n\in\textbf{N};\;E\left[X_{t_0+n}-X_{t_0}\right]^2
&>&k_2 C \|h_n\|^{2\alpha}-l_2\|h_n\|^{2H(t_0)}\\
&\geq& C'\|h_n\|^{2\alpha}
\end{eqnarray*}
hence, by lemma \ref{majloclem}
\[
\alpha\geq\alpha(t_0)
\]
and therefore
\begin{equation}\label{upponct2}
\alpha(t_0)\leq \beta(t_0)
\end{equation}
\end{itemize}
We can restate the upper bounds (\ref{upponct1}) and (\ref{upponct2})
of the pointwise H\"older exponent of $X$ at $t_0$
\begin{equation}
\alpha(t_0)\leq \beta(t_0)\wedge H(t_0)
\textrm{ (resp. }\beta(t_0)\wedge H_i(t_0)\textrm{ )}
\end{equation}

\subsubsection{Upper bound for the local H\"older exponent}
By (\ref{locinfpoint}), any upper bound for the pointwise H\"older exponent is an
upper bound for the local H\"older exponent. But we can improve on this result in
the case $\tilde{\beta}(t_0)<H(t_0)$.
We first give an analogous of lemma \ref{majloclem} for the local exponent
\begin{lemma}\label{majpointlem}
Let $X=\left\{X_t;t\in\textbf{R}_{+}^N\right\}$ be a Gaussian process.
Assume there exists $\mu\in (0,1)$ such that for all $\epsilon>0$, there exist
two sequences $\left(h_n\right)_{n\in\textbf{N}}$ and $\left(l_n\right)_{n\in\textbf{N}}$
of $\left(\textbf{R}_{+}^N\right)^{*}$ converging to 0, and a constant $c>0$
such that
\[
\forall n\in\textbf{N};\;
E\left[X_{t_0+h_n}-X_{t_0+l_n}\right]^2\geq c\|h_n-l_n\|^{2\mu+\epsilon}
\]
Then we have almost surely
\[
\tilde{\alpha}(t_0)\leq\mu
\]
\end{lemma}

\proof
Let $\epsilon>0$ and consider two sequences $(h_n)_{n\in\textbf{N}}$
and $(l_n)_{n\in\textbf{N}}$ as in the statement.\\
For all $n\in\textbf{N}$, the law of the random variable 
$\frac{X_{t_0+h_n}-X_{t_0+l_n}}{\|h_n-l_n\|^{\mu+\epsilon}}$ is
$\mathcal{N}(0,\sigma_n^2)$.\\
From the assumption, we have $\sigma_n\rightarrow +\infty$ 
as $n\rightarrow +\infty$.\\
Then, for all $\lambda>0$,
\begin{eqnarray*}
P\left\{\frac{\|h_n-l_n\|^{\mu+\epsilon}}{|X_{t_0+h_n}-X_{t_0+l_n}|}<\lambda
\right\}&=&P\left\{\frac{|X_{t_0+h_n}-X_{t_0+l_n}|}{\|h_n-l_n\|^{\mu+\epsilon}}
>\frac{1}{\lambda}\right\}\\
&=&\int_{|x|>\frac{1}{\lambda}}\frac{1}{\sqrt{2\pi}\sigma_n}
\exp\left(-\frac{x^2}{2\sigma_n^2}\right).dx\\
&=&\frac{1}{2\pi}\int_{|x|>\frac{1}{\lambda\sigma_n}}
\exp\left(-\frac{x^2}{2}\right).dx\stackrel{n\rightarrow +\infty}
{\longrightarrow} 1
\end{eqnarray*}
Therefore the sequence 
$\left(\frac{\|h_n-l_n\|^{\mu+\epsilon}}{|X_{t_0+h_n}-X_{t_0+l_n}|}\right)_{n\in\textbf{N}}$
converges to $0$ in probability. then there exists a subsequence which converges
to $0$ almost surely.
Then we have almost surely $\tilde{\alpha}(t_0)\leq\mu+\epsilon$. Taking
$\epsilon\in\textbf{Q}_{+}$, the result follows.

\fin

Let $\alpha\in(\tilde{\beta}(t_0);H(t_0))$ (resp.
$\alpha\in(\tilde{\beta}(t_0);H_i(t_0))$).
As
\begin{eqnarray*}
\limsup_{\rho\rightarrow 0}\sup_{s,t\in B(t_0,\rho)}
\frac{\left|H(t)-H(s)\right|}{\|t-s\|^{\alpha}}=+\infty
\end{eqnarray*}
for all $M>0$, there exists $\rho_0>0$ such that 
\begin{eqnarray*}
\forall \rho<\rho_0;
\exists s,t\in B(t_0,\rho);\; \left|H(t)-H(s)\right|>M\|t-s\|^{\alpha}
\end{eqnarray*}
Therefore we can construct two sequences $\left(h_n\right)$ and $\left(l_n\right)$
converging to $0$ such that
\begin{eqnarray*}
\forall n\in\textbf{N};\;\left|H(t_0+h_n)-H(t_0+l_n)\right|>M\|h_n-l_n\|^{\alpha}
\end{eqnarray*}
By lemma \ref{majpointlem}, we can deduce
\begin{equation}
\tilde{\alpha}(t_0)\leq\tilde{\beta}(t_0)
\end{equation}

\subsection{Directional H\"older exponents}
One may also define directional pointwise and local H\"older
exponents in the direction
$u\in\mathcal{U}=\left\{u\in\textbf{R}^N;\|u\|=1\right\}$ by
\begin{eqnarray*}
\alpha_u(t_0)&=&\sup\left\{\alpha; \lim_{\rho\rightarrow 0}\frac{|X_{t_0+\rho.u}-X_{t_0}|}{\rho^{\alpha}}=0\right\}
\end{eqnarray*}
and
\begin{eqnarray*}
\tilde{\alpha}_u(t_0)=\sup\left\{\alpha;\;
\limsup_{\rho\rightarrow 0}\sup_{\scriptstyle s,t\in B(t_0,\rho)\atop
\scriptstyle s,t\in t_0+\textbf{R}.u}\frac{|X_t-X_s|}{\|t-s\|^{\alpha}}<\infty\right\}
\end{eqnarray*}
As previously, for all $u\in\mathcal{U}$, we have
\begin{equation}
\tilde{\alpha}_u(t_0)\leq \alpha_u(t_0)
\end{equation}
Moreover, we can see easily that for all $u\in\mathcal{U}$, we have 
\begin{eqnarray}
\alpha(t_0)\leq\alpha_u(t_0)\textrm{ and }
\tilde{\alpha}(t_0)\leq\tilde{\alpha}_u(t_0)
\end{eqnarray}

\begin{proposition} \label{propexpdiriso}
Let $X=\left\{X_t; t\in\textbf{R}_{+}^N\right\}$ be a multifractional Brownian field.
For all $t_0\in\textbf{R}_{+}^N$ and all $u\in\mathcal{U}$, the local H\"older 
exponent of $X$ at $t_0$ in the direction $u$ is
almost surely given by
\begin{eqnarray}
\tilde{\alpha}_u(t_0)=\tilde{\beta}_u(t_0)\wedge H(t_0)
\end{eqnarray}
and the pointwise H\"older exponent of $X$ at $t_0$ in the direction $u$ 
satisfies almost surely
\begin{eqnarray}
\alpha_u(t_0)=\beta_u(t_0)\wedge H(t_0)
\end{eqnarray}
where $\beta_u(t_0)$ and $\tilde{\beta}_u(t_0)$ denote the pointwise and local
H\"older exponents of $H$ at $t_0$ in the direction $u$.
\end{proposition}

\proof
Let $t_0\in\textbf{R}_{+}^N$, $u\in\mathcal{U}$ and
consider the stochastic process
\[
\tilde{X}=\left\{\tilde{X}_{\rho}=X_{t_0+\rho.u};\rho>0\right\}
\]
By definition, $\alpha_u(t_0)$ and $\tilde{\alpha}_u(t_0)$ are respectively the
pointwise and local H\"older exponents of $\tilde{X}$ at $0$.\\
Let $\tilde{H}(\rho)=H(t_0+\rho.u)$. We have
\begin{eqnarray*}
\frac{1}{2}E\left[\tilde{X}_{\rho}-\tilde{X}_{\eta}\right]^2
&=&\frac{1}{2}E\left[X_{t_0+\rho.u}-X_{t_0+\eta.u}\right]^2\\
&=&D\left[2\tilde{H}(\rho)\right] |\rho-\eta|^{2\tilde{H}(\rho)}
+\frac{\partial^2\varphi}{\partial
x^2}\left(2\tilde{H}(\rho);\|t_0+\rho.u\|\right)\times
\left(\tilde{H}(\rho)-\tilde{H}(\eta)\right)^2\\
&&+o\left(\tilde{H}(\rho)-\tilde{H}(\eta)\right)^2
+o\left(|\rho-\eta|^{2\tilde{H}(\rho)}\right)
\end{eqnarray*}
by corollary \ref{corisoasymp}.\\
Then, using the same method as in proposition \ref{propexpiso}, the result follows.
\fin

\begin{proposition} \label{propexpdirdrap}
Let $X=\left\{X_t; t\in\textbf{R}_{+}^N\right\}$ be a multifractional Brownian
sheet.
For all $t_0\in\textbf{R}_{+}^N$, the local H\"older exponent of $X$ at $t_0$ in
the direction $\epsilon_i$ is
almost surely given by
\begin{eqnarray}
\tilde{\alpha}_{\epsilon_i}(t_0)=\tilde{\beta}_{\epsilon_i}(t_0)\wedge H_i(t_0)
\end{eqnarray}
and the pointwise H\"older exponent of $X$ at $t_0$ in the direction $u$ 
satisfies almost surely
\begin{eqnarray}
\alpha_{\epsilon_i}(t_0)=\beta_{\epsilon_i}(t_0)\wedge H_i(t_0)
\end{eqnarray}
where $\beta_{\epsilon_i}(t_0)$ and $\tilde{\beta}_{\epsilon_i}(t_0)$ denote the
pointwise and local H\"older exponents of $H$ at $t_0$ in the direction 
$\epsilon_i$.
\end{proposition}

\proof
As in the proof of lemma \ref{kolmdrap1}, there exists a constant $M>0$ and a
one-parameter mBm $X^{(i)}$ such that
\begin{eqnarray*}
E\left[X_{t_0+\rho.\epsilon_i}-X_{t_0+\eta.\epsilon_i}\right]^2
=M.E\left[X^{(i)}_{t_0^{(i)}+\rho.\epsilon_i}-X^{(i)}_{t_0^{(i)}+\eta.\epsilon_i}\right]^2
\end{eqnarray*}
Then, using the same method as in proposition \ref{propexpdrap}, the result follows.
\fin

\subsection{Application of Dudley's theory}
Another way to study the regularity of our processes is to examine the behavior 
around zero of the modulus of continuity 
\[
\omega_{X,T}(\delta)=\sup_{s,t\in T;\\ d(s,t)\leq\delta} |X_s-X_t|
\]

When the process studied is Gaussian, it is convenient to consider the
pseudo-metric 
\[
d(s,t)=E\left[X_s-X_t\right]^2
\]
As usually, we define the ball of radius $r>0$ about $t\in T$ by
\[
\mathcal{B}_d (t,r)=\left\{s\in T; d(s,t)<r\right\}
\]
and we say that $(T,d)$ is totally bounded if for all $\epsilon>0$,
there exists $t_1,\dots,t_m\in T$ such that 
\[
T\subset \bigcup_{i=1}^m \mathcal{B}_d (t_i,\epsilon)
\]
When $(T,d)$ is totally bounded, we can define the metric entropy 
$\epsilon\mapsto D(\epsilon, T,d)$ where $D(\epsilon, T,d)$ is minimum number
of balls of radius $\epsilon$ required to cover $T$.

The following theorem allows to improve on the results of the previous
paragraph. 

\begin{theorem}[Dudley's Theorem]
Consider a centered Gaussian process $X=\left\{X_t; t\in\textbf{R}_{+}^N\right\}$
indexed by the pseudo-metric space $(T,d)$.
If $(T,d)$ is totally bounded and if $\int_0^1 \sqrt{\ln D(r,T,d)} dr <\infty$,
then $X$ has a continuous modification $Y=\left\{Y_t; t\in\textbf{R}_{+}^N\right\}$.
Moreover, there exists a universal constant $C>0$ such that
\[
\limsup_{\delta\rightarrow 0^{+}} \frac{\omega_{Y,T}(\delta)}
{\int_0^{\delta} \sqrt{\ln D(\frac{r}{2},T,d)} dr
+C \delta\sqrt{\ln\ln\frac{1}{\delta}}}
\leq 24
\]

\end{theorem}
To apply this result, we first need to verify the assumptions about
the metric entropy.

\begin{lemma}
Let $T\subset\textbf{R}_{+}^N$ measurable and $d$ a pseudo-metric on $T$.
If there exists $C>0$ and $\alpha>0$ such that
\[
\forall s,t\in T;\; d(s,t)\leq C \|s-t\|^{\alpha}
\]
then there exists $r_0>0$ such that $\forall r\in [0,r_0]$,
\[
D(r,T,d)\leq C^{\frac{N}{\alpha}}.Leb(T).r^{-\frac{N}{\alpha}}
\]
\end{lemma}

We saw previously that the $2$ multi-parameter extensions of the mBm,
satisfy the assumption of this lemma with $T=[a,b]$.
Then there exists $C_{a,b}>0$, $\alpha=\alpha_{a,b}>0$ and $r_0>0$ 
such that for all $r\in[0,t_0]$,
\[
D(r,[a,b],d)\leq C^{\frac{N}{\alpha}}.Leb([a,b]).r^{-\frac{N}{\alpha}}
\]
As a consequence, $([a,b],d)$ is totally bounded and in the neighborhood of $0$,
we have
\[
\sqrt{\ln D(r,T,d)}\leq\sqrt{K-\frac{N}{\alpha}\ln r}
\]
therefore the integral $\int_0^1 \sqrt{\ln D(r,T,d)} dr$ is finite and we can apply 
Dudley's theorem.
We get

\begin{proposition}
Let $X=\left\{X_t; t\in\textbf{R}_{+}^N\right\}$ be one of the multi-parameter
extension of the mBm. For all $a\prec b$, there exists $C_{a,b}>0$ and
$\alpha=\alpha_{a,b}>0$ such that
\[
\limsup_{\epsilon\rightarrow 0}\frac{\sup_{\|s-t\|\leq\epsilon}|X_s-X_t|}
{\epsilon^{\alpha}\sqrt{\ln\frac{1}{\epsilon}}}
\leq C_{a,b}\sqrt{N}
\]
\end{proposition}

\proof
First of all, we study the quotient
\[
\frac{\int_0^{\delta} \sqrt{\ln D(\frac{r}{2},T,d)} dr
+C \delta\sqrt{\ln\ln\frac{1}{\delta}}}{\delta\sqrt{\ln\frac{1}{\delta}}}
=\frac{N(\delta)}{D(\delta)}
\]
The derivative of the numerator is
\begin{eqnarray*}
N'(\delta)&=&\sqrt{\ln D(\frac{\delta}{2},T,d)}+C\left[\sqrt{\ln\ln\frac{1}{\delta}}
-\frac{1}{2\ln\frac{1}{\delta}\sqrt{\ln\ln\frac{1}{\delta}}}\right]\\
&\leq&\sqrt{K-\frac{N}{\alpha}\ln\delta}+C\sqrt{\ln\ln\frac{1}{\delta}}
\left[1-\frac{1}{2\ln\frac{1}{\delta}\ln\ln\frac{1}{\delta}}\right]
\end{eqnarray*}
and the derivative of the denominator is
\begin{eqnarray*}
D'(\delta)&=&\sqrt{\ln\frac{1}{\delta}}\left[1-\frac{1}{2\ln\frac{1}{\delta}}\right]\\
&\sim&\sqrt{\ln\frac{1}{\delta}}
\end{eqnarray*}
Then we have
\[
\limsup_{\delta\rightarrow 0^{+}}\frac{N'(\delta)}{D'(\delta)}
\leq \sqrt{\frac{N}{\alpha}}
\]
and by a L'Hopital's rule type argument,
\[
\limsup_{\delta\rightarrow 0^{+}}\frac{N(\delta)}{D(\delta)}
\leq \sqrt{\frac{N}{\alpha}}
\]
Then we have
\[
\limsup_{\delta\rightarrow 0}\frac{\sup_{d(s,t)\leq\delta}|X_s-X_t|}
{\delta\sqrt{\ln\frac{1}{\delta}}}\leq 24\sqrt\frac{N}{\alpha}
\]
The problem is now to transform $\sup_{d(s,t)\leq\delta}$ into
$\sup_{\|s-t\|\leq\delta}$.\\
To do this, we write 
\[
\forall s,t\in[a,b];\; d(s,t)\leq C_{a,b}\|s-t\|^{\alpha}
\]
then $\|s-t\|\leq\epsilon=\frac{\delta^{\frac{1}{\alpha}}}{C_{a,b}^{\frac{1}{\alpha}}}$
implies $d(s,t)\leq\delta$
and we get
\[
\limsup_{\delta\rightarrow 0}\frac{\sup_{\|s-t\|\leq\epsilon}|X_s-X_t|}
{C_{a,b}\epsilon^{\alpha}\sqrt{\alpha\ln\frac{1}{\epsilon}}}
\leq 24\sqrt\frac{N}{\alpha}
\]
which gives the expected result.
\fin

This result is more powerful than knowledge of H\"older exponents.
It gives the behavior of $\left|X_s-X_t\right|$ in a ball around $t_0$.

\section{Locally asymptotic self-similarity}
Extending fBm into multifractional Brownian motion implies the loss of the two 
properties of self-similarity and stationarity of increments. However, a weak
form of self-similarity remains, called locally asymptotic self-similarity (see
\cite{ayach}, \cite{benassi}).
As we will see, this property still holds for the two kinds of extension of mBm
in $\textbf{R}^N$.

\begin{theorem}\label{thselfiso}
Let $X=\left\{X_t;t\in\textbf{R}_{+}^N\right\}$ be a multifractional
Brownian field.\\
For all $t_0\in\textbf{R}_{+}^N$, the law of the process 
$Y^{\alpha}(\rho)=\left\{Y^{\alpha}_u(\rho)=
\frac{X_{t_0+\rho u}-X_{t_0}}{\rho^{\alpha}}; u\in\textbf{R}_{+}^N\right\}$ 
converges weakly if one of the following two conditions holds
\begin{enumerate}
\item $\alpha=H(t_0)$ and $H(t_0)<\inf_{u,v}\beta_{uv}(t_0)$\\
where $\beta_{uv}(t_0)=\sup\left\{\alpha;\;\lim_{\rho\rightarrow 0}
\frac{|H(t_0+\rho u)-H(t_0+\rho v)|}{\rho^{\alpha}}=0\right\}$.\\
Then, the limit measure is the law of a fractional Brownian field 
with parameter $H(t_0)$.
\item $\alpha=\inf_{u,v}\beta_{uv}(t_0)$,
$H(t_0)>\inf_{u,v}\beta_{uv}(t_0)$
and for all $u,v\in\textbf{R}_{+}^N$, the following limit exists
\[
\lim_{\rho\rightarrow 0}\frac{|H(t_0+\rho u)-H(t_0+\rho v)|}
{\rho^{\inf_{u,v}\beta_{uv}(t_0)}}=\Gamma(u,v)
\]
with $(u,v)\mapsto\frac{\Gamma(u,v)}{\|u-v\|^{2\beta}}$ bounded on $[a,b]^2$ for
some $\beta>0$.\\
The limit measure is the law of a Gaussian process 
$Y^{\inf_{u,v}\beta_{uv}(t_0)}$ such that
\[
E\left[Y^{\inf_{u,v}\beta_{uv}(t_0)}_u-Y^{\inf_{u,v}\beta_{uv}(t_0)}_v\right]^2
=K_{t_0}\left[\Gamma(u,v)\right]^2
\]
\end{enumerate}
\end{theorem}

\begin{remark}
As in the Levy fBm's case in proposition \ref{deltaiso}, the same result as
theorem \ref{thselfiso} can be stated for the increments $\Delta X$ defined in
section \ref{increments}.
The law of the process
$Y^{\alpha}(\rho)=\left\{Y^{\alpha}_u(\rho)=
\frac{\Delta X_{t_0,t_0+\rho u}}{\rho^{\alpha}}; u\in\textbf{R}_{+}^N\right\}$
converges weakly under the same assumptions. 
\end{remark}

In the case $N=1$, for all $u,v\in\textbf{R}_{+}$, we have
$\beta_{uv}(t_0)=\beta(t_0)$. Therefore, theorem \ref{thselfiso} has a simpler
statement. The two cases to be considered, depend of the comparison between
$H(t_0)$ and the pointwise exponent $\beta(t_0)$ of $H$.
 
The following example shows that the limit considered in the second case, can be
non trivial.
\begin{example}
In the case $N=1$, let $H(t)=\frac{3}{4}+t^{\frac{1}{2}}$ for $t\in[0,\frac{1}{4}]$.\\
For $t_0=0$, we compute, for all $u,v$ and $\rho>0$
\[
\frac{|H(\rho.u)-H(\rho.v)|}{\rho^{\frac{1}{2}}}=
|u^{\frac{1}{2}}-v^{\frac{1}{2}}|<|u-v|^{\frac{1}{2}}
\]
The limit measure is the law of a centered Gaussian process $Y$ such that
\[
E\left[Y_u-Y_v\right]^2=K_0\;\left(u^{\frac{1}{2}}-v^{\frac{1}{2}}\right)^2
\]
ie
\[
E\left[Y_u Y_v\right]=K_0\;u^{\frac{1}{2}}v^{\frac{1}{2}}
\]
\end{example}

\begin{theorem} \label{thasympdrap}
Let $X=\left\{X_t;t\in\textbf{R}_{+}^N\right\}$ be a multifractional
Brownian sheet.\\
The law of the process $Y^{\alpha}(\rho)=\left\{Y^{\alpha}_u(\rho)=\frac{\Delta
X_{t_0,t_0+\rho u}}{\rho^{\sum_i \alpha_i}}; u\in\textbf{R}_{+}^N\right\}$ 
converges weakly if for all $i\in\left\{1,\dots,N\right\}$, one of the following
two conditions holds
\begin{enumerate}
\item $\alpha_i=H_i(t_0)$ and $H_i(t_0)<\inf_{u,v}\beta^i_{uv}(t_0)$\\
where $\beta^i_{uv}(t_0)=\sup\left\{\alpha;\;\lim_{\rho\rightarrow 0}
\frac{|H_i(t_0+\rho u)-H_i(t_0+\rho v)|}{\rho^{\alpha}}=0\right\}$.
\item $\alpha_i=\inf_{u,v}\beta^i_{uv}(t_0)$,
$H_i(t_0)>\inf_{u,v}\beta^i_{uv}(t_0)$
and
\[
\lim_{\rho\rightarrow 0}\frac{|H_i(t_0+\rho u)-H_i(t_0+\rho v)|}
{\rho^{\inf_{u,v}\beta^i_{uv}(t_0)}}=\Gamma_i(u,v)
\]
with $(u,v)\mapsto\frac{\Gamma_i(u,v)}{\|u-v\|^{2\beta_i}}$ bounded on $[a,b]^2$
for some $\beta_i>0$.
\end{enumerate}
\end{theorem}

As usually, the proof of weak convergence proceeds in two steps. First, we need 
to show finite dimensional convergence, and then,
use a tightness argument.
Lemma 14.2 and theorem 14.3 in \cite{kal}, for instance, allow then to conclude.

\subsection{Finite dimensional convergence}
As the considered processes are Gaussian, we only have to show the convergence
of covariance functions.

\subsubsection{Multifractional Brownian field}
By (\ref{covisoasymp}), we compute
\begin{eqnarray}\label{fddiso}
\rho^{2\alpha}E\left[Y^{\alpha}_u(\rho)-Y^{\alpha}_v(\rho)\right]^2
&=&E\left[X_{t_0+\rho u}-X_{t_0+\rho v}\right]^2\nonumber\\
&=&D\left[H(t_0+\rho u)+H(t_0+\rho v)\right]\times\|\rho.(u-v)\|^{H(t_0+\rho
u)+H(t_0+\rho v)}\nonumber\\
&&+\frac{\partial^2\varphi}{\partial x^2}\left(2H(t_0+\rho u);\|t_0+\rho u\|\right)
\times\left(H(t_0+\rho u)-H(t_0+\rho v)\right)^2\nonumber\\
&&+o\left(\|\rho.(u-v)\|^2\right)+o\left(H(t_0+\rho u)-H(t_0+\rho v)\right)^2
\end{eqnarray}
To show that 
\[
\rho^{H(t_0+\rho u)+H(t_0+\rho v)}\sim\rho^{2H(t_0)}
\]
in the neighborhood of $\rho=0$, we study
\begin{eqnarray*}
\left[H(t_0+\rho u)+H(t_0+\rho v)-2H(t_0)\right]\;\ln\rho&=&
\frac{H(t_0+\rho u)-H(t_0)}{\|\rho.u\|^{\alpha}}\times\|\rho.u\|^{\alpha}\ln\rho\\
&&+\frac{H(t_0+\rho v)-H(t_0)}{\|\rho.v\|^{\alpha}}\times\|\rho.v\|^{\alpha}\ln\rho
\end{eqnarray*}
for $\alpha<\beta(t_0)$.\\
As $\left(u;\rho\right)\mapsto\|\rho.u\|^{\alpha}\ln\rho$ is bounded on
$[a,b]\times[0,1]$
and 
\[
\forall u\in[a,b];\; \frac{H(t_0+\rho u)-H(t_0)}{\|\rho.u\|^{\alpha}}
\stackrel{\rho\rightarrow 0}{\longrightarrow}0
\]
we have
\[
\left[H(t_0+\rho u)+H(t_0+\rho v)-2H(t_0)\right]\;\ln\rho
\stackrel{\rho\rightarrow 0}{\longrightarrow}0
\]
Therefore, in the neighborhood of $\rho=0$, the first term of (\ref{fddiso}) is equivalent to
\[
D\left[2H(t_0)\right]\|u-v\|^{2H(t_0)}\times\rho^{2H(t_0)}
\]
and the second to
\[
\frac{\partial^2\varphi}{\partial x^2}\left(2H(t_0);\|t_0\|\right)
\times\left(H(t_0+\rho u)-H(t_0+\rho v)\right)^2
\]
Let 
$\beta_{uv}(t_0)=\sup\left\{\alpha;\lim_{\rho\rightarrow 0}\frac{|H(t_0+\rho
u)-H(t_0+\rho v)|}{\rho^{\alpha}}=0\right\}$.
We have to distinguish the two following cases
\begin{itemize}
\item if $H(t_0)<\inf_{u,v}\beta_{uv}(t_0)$, by definition of
$\beta_{uv}(t_0)$,
\begin{eqnarray*}
\forall u,v\in\textbf{R}_{+}^N;\;
\lim_{\rho\rightarrow 0}\frac{\left|H(t_0+\rho u)-H(t_0+\rho v)\right|}
{\rho^{H(t_0)}}=0
\end{eqnarray*}
Therefore
\begin{eqnarray*}
\forall u,v\in\textbf{R}_{+}^N;\;
E\left[Y^{H(t_0)}_u(\rho)-Y^{H(t_0)}_v(\rho)\right]^2
\stackrel{\rho\rightarrow 0}{\longrightarrow}
\underbrace{D\left[2H(t_0)\right]\|u-v\|^{2H(t_0)}}
_{E\left[B^{H(t_0)}_u-B^{H(t_0)}_v\right]^2}
\end{eqnarray*}
where $B^{H(t_0)}$ denotes fractional Brownian field of parameter $H(t_0)$.
\item if $H(t_0)>\inf_{u,v}\beta_{uv}(t_0)$,\\
for all $\alpha<\inf_{u,v}\beta_{uv}(t_0)$, as
\begin{eqnarray*}
\forall u,v\in\textbf{R}_{+}^N;\;
\lim_{\rho\rightarrow 0}\frac{\left|H(t_0+\rho u)-H(t_0+\rho v)\right|}
{\rho^{\alpha}}=0
\end{eqnarray*}
we have
\begin{eqnarray*}
\forall u,v\in\textbf{R}_{+}^N;\;
\frac{1}{\rho^{2\alpha}}E\left[X_{t_0+\rho u}-X_{t_0+\rho v}\right]^2
\stackrel{\rho\rightarrow 0}{\longrightarrow} 0
\end{eqnarray*}
Moreover, since there exists $u,v\in\textbf{R}_{+}^N$ 
such that $H(t_0)>\beta_{uv}(t_0)$, we can
consider $\alpha\in\left(\beta_{uv}(t_0);H(t_0)\right)$.
The limit
\[
\limsup_{\rho\rightarrow 0}\frac{\left|H(t_0+\rho u)-H(t_0+\rho v)\right|}
{\rho^{\alpha}}=+\infty
\]
implies
\begin{eqnarray*}
\exists u,v\in\textbf{R}_{+}^N;\;
\limsup_{\rho\rightarrow 0}\frac{1}{\rho^{2\alpha}}E\left[X_{t_0+\rho u}-X_{t_0+\rho v}\right]^2
=+\infty
\end{eqnarray*}
Therefore $E\left[Y^{\alpha}_u(\rho)-Y^{\alpha}_v(\rho)\right]^2$ admits a
limit for all $u,v\in\textbf{R}_{+}^N$ when
$\rho\rightarrow 0$ {\bf if and only if} 
\begin{eqnarray*}
\left\{\begin{array}{ll}
\alpha=\inf_{u,v}\beta_{uv}(t_0)\\
\textrm{and}\\ 
\lim_{\rho\rightarrow 0}\frac{|H(t_0+\rho u)-H(t_0+\rho v)|}
{\rho^{\inf_{u,v}\beta_{uv}(t_0)}}=\Gamma(u,v)\in\textbf{R}_{+}^{*}
\end{array}\right.
\end{eqnarray*}
In that case, we have for all $u$, $v$ in $\textbf{R}^N_{+}$,
\begin{eqnarray*}
E\left[Y_u^{\alpha}(\rho)-Y_v^{\alpha}(\rho)\right]^2
\stackrel{\rho\rightarrow 0}{\longrightarrow}
\frac{\partial^2\varphi}{\partial x^2}\left(2H(t_0);\|t_0\|\right)
\left[\Gamma(u,v)\right]^2
\end{eqnarray*}

\begin{remark}
We can see easily that
\begin{equation}
\beta_{\frac{u}{\|u\|}}(t_0)\wedge\beta_{\frac{v}{\|v\|}}(t_0)
\leq\beta_{uv}(t_0)
\end{equation}
hence
\begin{equation}
\inf_{u\in\mathcal{U}}\beta_u(t_0)\leq\inf_{u,v}\beta_{uv}(t_0)
\end{equation}

Conversely, assume there exist $u,v\in\mathcal{U}$ such that 
$\beta_u(t_0)<\beta_v(t_0)$, and let $\alpha\in(\beta_u(t_0);\beta_v(t_0))$.
The inequality
\[
\frac{|H(t_0+\rho u)-H(t_0)|}{\rho^{\alpha}}\leq
\frac{|H(t_0+\rho u)-H(t_0+\rho v)|}{\rho^{\alpha}}
+\frac{|H(t_0+\rho v)-H(t_0)|}{\rho^{\alpha}}
\]
implies
\[
\limsup_{\rho\rightarrow 0}\frac{|H(t_0+\rho u)-H(t_0+\rho v)|}{\rho^{\alpha}}
=+\infty
\]
and therefore $\alpha>\beta_{uv}(t_0)$.
Then $\inf_{u,v}\beta_{uv}(t_0)\leq\inf_{u\in\mathcal{U}}\beta_u(t_0)$,
which gives
\begin{equation}
\inf_{u,v}\beta_{uv}(t_0)=\inf_{u\in\mathcal{U}}\beta_u(t_0)
\end{equation}
\end{remark}

\end{itemize}

\subsubsection{Multifractional Brownian sheet}
In the non-isotropic case, using remark \ref{prodtens},
consider $N$ independent mBm $X^{(i)}$ with parameter function $H_i$
\begin{eqnarray*}
E\left[\Delta X_{t_0,t_0+\rho u} \Delta X_{t_0,t_0+\rho v}\right]
=\prod_{i=1}^N E\left[\left(X^{(i)}_{t^{(i)}_0+\rho u^{(i)}}
-X^{(i)}_{t^{(i)}_0}\right)
\left(X^{(i)}_{t^{(i)}_0+\rho v^{(i)}}-X^{(i)}_{t^{(i)}_0}\right)\right]
\end{eqnarray*}
As in the isotropic case, for all $i\in\left\{1,\dots,N\right\}$, consider
\[
\beta^i_{uv}(t_0)=\sup\left\{\alpha;\;\lim_{\rho\rightarrow 0}
\frac{|H_i(t_0+\rho u)-H_i(t_0+\rho v)|}{\rho^{\alpha}}=0\right\}
\]
Each process $X^{(i)}$ is locally asymptoticly self-similar, therefore
\[
E\left[\frac{X^{(i)}_{t^{(i)}_0+\rho u^{(i)}}
-X^{(i)}_{t^{(i)}_0}}{\rho^{\alpha_i}}
\times\frac{X^{(i)}_{t^{(i)}_0+\rho v^{(i)}}
-X^{(i)}_{t^{(i)}_0}}{\rho^{\alpha_i}}\right]
\stackrel{\rho\rightarrow 0}{\longrightarrow}
E\left[Y^{\alpha_i}_{u}Y^{\alpha_i}_{v}\right]
\]
where $Y^{\alpha_i}$ denotes 
\begin{itemize}
\item fractional Brownian motion of parameter $\alpha_i=H_i(t_0)$, in the case
$H_i(t_0)<\inf_{u,v}\beta^i_{uv}(t_0)$,
\item the centered Gaussian process such that
\[
E\left[Y^{\alpha_i}_u-Y^{\alpha_i}_v\right]^2
=K_{t_0}\left[\Gamma_i(u,v)\right]^2
\]
where $\alpha_i=\inf_{u,v}\beta^i_{uv}(t_0)$,
in the case $H_i(t_0)>\inf_{u,v}\beta^i_{uv}(t_0)$ and
\[
\lim_{\rho\rightarrow 0}\frac{|H_i(t_0+\rho u)-H_i(t_0+\rho v)|}
{\rho^{\inf_{u,v}\beta^i_{uv}(t_0)}}=\Gamma_i(u,v)
\]
with $\Gamma_i$ bounded on $[a,b]^2$.
\end{itemize}
Then we conclude
\[
E\left[\frac{\Delta X^{(i)}_{t^{(i)}_0,t^{(i)}_0+\rho u^{(i)}}}
{\rho^{\sum_i \alpha_i}}
\times\frac{\Delta X^{(i)}_{t^{(i)}_0,t^{(i)}_0+\rho v^{(i)}}}
{\rho^{\sum_i \alpha_i}}\right]
\stackrel{\rho\rightarrow 0}{\longrightarrow}
E\left[Y^{\alpha}_{u}Y^{\alpha}_{v}\right]
\]
where $Y^{\alpha}=\prod_{i=1}^N Y^{\alpha_i}$.

\subsection{Tightness of laws}
The study of weak convergence is well-known for stochastic processes indexed
by $\textbf{R}_{+}$. A comprehensive review was made by Billingsley 
(cf \cite{bill}) 
for a compact set of index ($[0,1]$). In (\cite{ks}), Karatzas and Shreeve stated
the same kind of results for the whole  $\textbf{R}_{+}$.
The case of $\textbf{R}_{+}^N$ can be found in (\cite{kal}) whose corollary 14.9
provides

\begin{proposition}\label{tightprop}
Consider a sequence of continuous processes 
$\left(X^{(n)}\right)_{n\in\textbf{N}}$ with $X^{(n)}=\left\{X^{(n)}_t;t\in\textbf{R}_{+}^N\right\}$
on $\left(\Omega,\mathcal{F},P\right)$ such that
\begin{enumerate}
\item there exists a positive constant $\nu$ such that
\[
\sup_{n\geq 1} E\left|X^{(n)}_0\right|^{\nu} <\infty
\]
\item for all $T>0$ and all $s$, $t$ in $[0,T]^N$, there exist
positive constants $\alpha$, $\beta$  and $C_T$ such that
\[
\sup_{n\geq 1} E\left|X^{(n)}_t-X^{(n)}_s\right|^{\alpha} \leq
C_T \|t-s\|^{N+\beta}
\]
\end{enumerate}
Then the probability measures $P_n\stackrel{\Delta}{=}P.\left(X^{(n)}\right)^{-1}$
on $\left(C\left(\textbf{R}_{+}^N\right),\mathcal{B}\left(C\left(\textbf{R}_{+}^N\right)\right)\right)$
form a tight sequence.
\end{proposition}

We verify the conditions of proposition \ref{tightprop}, in the case of mBm, in
the following sections.

\subsubsection{Multifractional Brownian field}
By (\ref{majcoviso}), there exist positive
constants $K_T$ and $L_T$ such that for all $u$, $v$ in $[0,T]^N$
\begin{eqnarray*}
\rho^{2\alpha}E\left[Y^{\alpha}_u(\rho)-Y^{\alpha}_v(\rho)\right]^2
&=&E\left[X_{t_0+\rho u}-X_{t_0+\rho v}\right]^2\\
&\leq& K_T\; \|\rho.(u-v)\|^{2H(t_0+\rho u)}\\
&&+L_T \left|H(t_0+\rho u)-H(t_0+\rho v)\right|^2
\end{eqnarray*}
Therefore,
\begin{eqnarray*}
E\left[Y^{\alpha}_u(\rho)-Y^{\alpha}_v(\rho)\right]^2
\leq K'_T\; \rho^{2(H(t_0)-\alpha)}.\|(u-v)\|^{2H(t_0)}
+L_T \frac{\left|H(t_0+\rho u)-H(t_0+\rho v)\right|^2}{\rho^{2\alpha}}
\end{eqnarray*}
\begin{itemize}
\item In the case $H(t_0)<\inf_{u,v}\beta_{uv}(t_0)$, there exists $M_T>0$ such that
\begin{eqnarray*}
E\left[Y^{H(t_0)}_u(\rho)-Y^{H(t_0)}_v(\rho)\right]^2
\leq M_T\; \|u-v\|^{2H(t_0)}
\end{eqnarray*}
\item In the case $H(t_0)>\inf_{u,v}\beta_{uv}(t_0)$, under the assumption 
\[
\lim_{\rho\rightarrow 0}\frac{\left|H(t_0+\rho u)-H(t_0+\rho v)\right|}
{\rho^{\inf_{u,v}\beta_{uv}(t_0)}}
=\Gamma(u,v)
\]
with $(u,v)\mapsto\frac{\Gamma(u,v)}{\|u-v\|^{2\beta}}$ bounded on $[a,b]^2$,
there exists $M_T>0$ such that 
\begin{eqnarray*}
E\left[Y^{\inf_{u,v}\beta_{uv}(t_0)}_u(\rho)-Y^{\inf_{u,v}\beta_{uv}(t_0)}_v(\rho)\right]^2
\leq M_T\; \|u-v\|^{2(\beta\wedge H(t_0))}
\end{eqnarray*}
\end{itemize}

Since the process $Y^{\alpha}$ is Gaussian, 
we get an exponent greater than $N$ in the usual way. 
Then we can conclude by proposition \ref{tightprop} that the laws of
$Y^{\alpha}$ are tight.

\subsubsection{Multifractional Brownian sheet}
In the same way as in paragraph \ref{parkolmdrap}
\begin{eqnarray*}
E\left[Y^{(\rho)}_u-Y^{(\rho)}_v\right]^2&=&\frac{1}{\rho^{2\sum_i \alpha_i}}
E\left[\prod_{i=1}^N \left(X^{(i)}_{t^{(i)}_0+\rho u^{(i)}}
-X^{(i)}_{t^{(i)}_0}\right)-\prod_{i=1}^N \left(X^{(i)}_{t^{(i)}_0+\rho v^{(i)}}
-X^{(i)}_{t^{(i)}_0}\right)\right]\\
&\leq&K \sum_i E\left[\frac{X^{(i)}_{t^{(i)}_0+\rho u^{(i)}}-X^{(i)}_{t^{(i)}_0+\rho
v^{(i)}}}{\rho^{\alpha_i}}\right]^2
\end{eqnarray*}
then, under the assumptions of theorem \ref{thasympdrap}, there exists a
positive constant $M_T$ such that
\begin{eqnarray*}
E\left[Y^{(\rho)}_u-Y^{(\rho)}_v\right]^2\leq M_T\;\|u-v\|^{2\min_i \alpha_i}
\end{eqnarray*}
We conclude as in the isotropic case.

\section*{Acknowledgement}
The author thanks Jacques L\'evy-V\'ehel for all their fruitful
discussions, especially about the H\"older regularity.

\bibliographystyle{plain}
\bibliography{style}

\end{document}